\documentclass[11pt, a4paper, reqno]{amsart}
\usepackage[english]{babel}
\usepackage{amsthm, fullpage}
\usepackage{graphicx}
\usepackage{amsfonts}
\usepackage{amsopn}
\usepackage{amsmath, amsfonts, amssymb, amscd, amsthm}
\usepackage{enumitem}
\usepackage{yhmath}
\usepackage{setspace}
\usepackage{mathtools}
\usepackage{mathrsfs}

\usepackage[all]{xy}

\newtheorem{theorem}{Theorem}[section]
\newtheorem{lemma}[theorem]{Lemma}
\newtheorem{proposition}[theorem]{Proposition}
\newtheorem{corollary}[theorem]{Corollary}
\theoremstyle{definition}
\newtheorem{definition}[theorem]{Definition}
\theoremstyle{remark}
\newtheorem{remark}[theorem]{Remark}
\theoremstyle{example}

\theoremstyle{note}

\numberwithin{equation}{section}

\DeclareMathOperator{\Hom}{Hom}

\DeclareMathOperator{\GL}{GL}
\DeclareMathOperator{\SL}{SL}

\DeclareMathOperator{\M}{M}

\DeclareMathOperator{\Aut}{Aut}
\DeclareMathOperator{\tr}{Tr}

\DeclareMathOperator{\coh}{H}
\DeclareMathOperator{\Rad}{Rad}
\DeclareMathOperator{\Char}{char}
\DeclareMathOperator{\Rep}{Rep}
\DeclareMathOperator{\End}{End}

\usepackage[colorlinks]{hyperref}

\usepackage{comment}

\usepackage{mathabx} 

\begin{document}

\title{ Dualizing involutions on Kazhdan-Patterson covers of $\GL(r)$}

\author{Ila Ahmad} 
\address{Ila Ahmad\\
Department of Mathematics\\
IISER Bhopal\\
Bhopal, Madhya Pradesh 462066, India}
\email{ila23@iiserb.ac.in}

\author{Kumar Balasubramanian} 

\address{Kumar Balasubramanian\\
Department of Mathematics\\
IISER Bhopal\\
Bhopal, Madhya Pradesh 462066, India}
\email{bkumar@iiserb.ac.in}

\author{Sanjeev Kumar Pandey} 
\address{Sanjeev Kumar Pandey\\
Department of Mathematics\\
IISER Bhopal\\
Bhopal, Madhya Pradesh 462066, India}
\email{sanjeevpandey@iiserb.ac.in}

\author{Ekta Tiwari} 
\address{Ekta Tiwari\\
Department of Mathematics\\
IISER Pune\\
Pune, Maharashtra 411008, India}
\email{ekta.tiwari@iiserpune.ac.in}

\author{Varsha Vasudevan} 
\address{Varsha Vasudevan\\
Department of Mathematics\\
IISER Bhopal\\
Bhopal, Madhya Pradesh 462066, India}
\email{varsha23@iiserb.ac.in}


\maketitle

\begin{abstract}
Let $F$ be a non-Archimedean local field of characteristic zero and $G=\GL(r,F).$ Let $r,n \geq 2$ and $\widetilde{G}$ be an $n$-fold Kazhdan-Patterson cover of $G$. Let $\pi$ be an irreducible smooth representation of $G$ and $\pi^{\vee}$ be the contragredient of $\pi$. Let $\tau$ be an involutive anti-automorphism of $G$ satisfying $\pi^{\tau}\simeq \pi^{\vee}$. In this case, we say that $\tau$ is a dualizing or duality involution. When $r=2$, it is known that any lift of the standard involution on $G$ to $\widetilde{G}$ is a dualizing involution if and only if $n=2$. In  this paper, we study the dualizing property of the lifts of the standard involution on $G$ to $\widetilde{G}$ and prove the following: (i) When $r \geq 2$ and $n \geq 3$,  we show that any lift of the standard involution to $\widetilde{G}$ is never a dualizing involution. (ii) When $r=2t$ and $n=2$, we show that any lift of the standard involution to $\widetilde{G}$ is a dualizing  involution. (iii) When $r=2t+1$ and $n=2$, we show that there is a unique lift of the standard involution to $\widetilde{G}$ which is a dualizing involution. 
\end{abstract}

\section{Introduction}

Let $F$ be a non-Archimedean local field of characteristic zero. Let $r\geq 2$ be a positive integer and $G=\GL(r,F)$. For $g\in G$, we let $^{\top}g$ denote the transpose of the matrix $g$, and $w_{0}$ to be the matrix with anti-diagonal entries equal to one. Let $\tau: G\rightarrow G$ be the map $\tau(g)=w_{0} {}^{\top}gw_{0}$. Clearly, $\tau$ is an anti-automorphism of $G$ and satisfies $\tau^{2}=1$. We call $\tau$ the \textit{standard involution} on $G$. Let $(\pi,V)$ be an irreducible smooth complex representation of $G$ and $(\pi^{\vee}, V^{\vee})$ be the smooth dual or the contragredient of $(\pi,V)$. For $\theta$ an anti-automorphism of $G$, we let $\pi^{\theta}$ to be the twisted representation of $G$ defined by \[\pi^{\theta}(g)=\pi(\theta(g^{-1})).\]
A well known result of Gelfand and Kazhdan is the following: 

\begin{theorem}[Gelfand-Kazhdan] Let $\tau$ be the standard involution on $G$. Then \[\pi^{\tau}\simeq \pi^{\vee}.\]
\end{theorem}
We refer the reader to Theorem 2 in \cite{GelKaj} for a proof of the above result. \\

If $\theta$ is any anti-automorphism of $G$ such that $\theta^{2}=1$, and satisfies $\pi^{\theta}\simeq \pi^{\vee}$, then we call $\theta$  a \textit{dualizing involution} or a \textit{duality involution}. The above result of Gelfand and Kazhdan says that the standard involution $\tau$ on $G$ is a dualizing involution. \\

For $ n\geq 2$, let $\widetilde{G}$ be a certain $n$-fold Kazhdan-Patterson cover of $G$ (see Section~\ref{Kazhdan-Patterson covers of GL(r)} below for a quick introduction or Chapter 0 in \cite{KazPat[2]} for a more detailed description of these covers). It is  known that (see Proposition 3.1 in \cite{KazPat}) the standard involution $\tau$ on $G$ has at least one lift to $\widetilde{G}$. A natural and interesting question that arises is the following: whether the lifts of the standard involution to $\widetilde{G}$ are also dualizing involutions. In the case when $\widetilde{G}$ is the $2$-fold cover of $G = \GL(2,F)$, it is known that any lift of $\tau$ to $\widetilde{G}$ is dualizing (see \cite{KumAji} and \cite{KumEkt}). In a more recent work \cite{Kumar-Renu-Sanjeev-Varsha}, it is proved that any lift of the standard involution $\tau$ of $G=\GL(2,F)$ to the $n$-fold cover $\widetilde{G}$ is a duality involution if and only if $n=2$. In other words, for $G=\GL(2,F)$ and $n\geq 3$, none of the lifts of $\tau$ to the $n$-fold cover $\widetilde{G}$ satisfy the dualizing property. These results for the $n$-fold covers of $\GL(2,F)$ motivate the analogous questions considered in this paper for the $n$-fold Kazhdan-Patterson covers of $G=\GL(r,F)$ for $n,r\geq 2$. \\

Before we state our results, we set up some notation. Let $F$ be a non-Archimedean local field of characteristic zero and $p$ be the residual degree of $F$. Let $r,n\geq 2$ and $G=\GL(r,F)$. Let $\mu_{n}$ be the group of $n$-th roots of unity in $F$. We assume that $|\mu_{n}|=n$. For $m\in \mathbb{Z}/n\mathbb{Z}$ and a certain $2$-cocycle $c^{(m)}\in \coh^{2}(G, \mu_{n})$, we let $\widetilde{G}^{(m)}$ denote the $n$-fold Kazhdan-Patterson cover of $G$ with respect to the $2$-cocycle $c^{(m)}$. When $m=0$, we simply write $\widetilde{G}^{(0)}=\widetilde{G}$ and $c^{(0)}=c$. We prove three main results about dualizing lifts of the standard involution to the cover $(\widetilde{G}, c)$. We state them below. 

\begin{theorem}\label{not duality}
Let $n \geq 3$ and $G = \GL(r,F)$ for $r\geq 2$. Let $\widetilde{G}$ be the $n$-fold Kazhdan-Patterson cover of $G$ with respect to the cocycle $c$. Let $\sigma$ be any lift of $\tau$ to $\widetilde{G}$. Then $\sigma$ is not a dualizing involution. That is, for any irreducible admissible genuine representation $\pi$ of $\widetilde{G}$, we have 
\[\pi^{\sigma}\not\simeq\pi^{\vee}.\]
\end{theorem}

\begin{theorem}\label{duality for even rank case}
Let $r =2t$ be an even integer and $G = \GL(r,F)$. Let $\widetilde{G}$ be the $2$-fold Kazhdan-Patterson cover of $G$ with respect to the cocycle $c$. Let $\sigma$ be any lift of $\tau$ to $\widetilde{G}$. Then $\sigma$ is a dualizing involution. That is, for any irreducible admissible genuine representation $\pi$ of $\widetilde{G}$, we have
\[\pi^{\sigma}\simeq \pi^{\vee}. \]
\end{theorem}

\begin{theorem}\label{unique duality for odd rank case}
Let $r =2t+1$ be an odd integer and $G = \GL(r,F)$. Let $\widetilde{G}$ be the $2$-fold Kazhdan-Patterson cover of $G$ with respect to the cocycle $c$. Then there exists a unique lift $\sigma$  of $\tau$ to $\widetilde{G}$ such that $\sigma$ is a dualizing involution. That is, for any irreducible admissible genuine representation $\pi$ of $\widetilde{G}$, we have 
\[\pi^{\sigma}\simeq \pi^{\vee}\]
for a unique lift $\sigma$ of $\tau$. 
\end{theorem}

Let $\theta$ be the anti-automorphism of $G$ given by $g\mapsto ~^{\top}g$. Let $\rho$ be any lift of $\theta$ to $\widetilde{G}$. 
We let $\sigma$ to be the lift of $\tau$ to $\widetilde{G}$ defined by $\tilde{g}\mapsto \tilde{w}_{0} \rho(\tilde{g}) \tilde{w}_{0}^{-1},$ where $\tilde{w}_{0}=(w_{0},1)\in \widetilde{G}$. It is easy to see that $\sigma$ is dualizing or non-dualizing if and only if the corresponding property holds for $\rho$. Under the assumption that $r, n\geq 3$ and $\gcd(n,p)=1$, we show that $\rho$ is not a duality involution by constructing $f\in C_{c}^{\infty}(\widetilde{G})$ such that the distribution characters of $\pi^{\rho}$ and $\pi^{\vee}$ evaluated at $f$ are distinct. This proof does not work in the case when $r=2$. For $r\geq 2$ and $n\geq 3$, we discuss another proof of the non-dualizing property of $\rho$ without any restrictions on the residual degree of $F$. On the other hand, to prove that $\rho$ is dualizing, the key idea is to show that the lift $\rho$ of $\theta$ satisfies a certain conjugacy property when restricted to the set of regular semisimple elements in $\widetilde{G}$. Adapting the results proved for the cover $(\widetilde{G},c)$ for the other covers $(\widetilde{G}^{(m)}, c^{(m)})$ follows \textit{mutatis mutandis} after some suitable adjustments (see Section~\ref{Twisted covers} for the precise statements in this setup). \\

The study of covering groups is not only interesting from a
representation-theoretic point of view, but also has important connections
with number theory and the Langlands program. In the past, having a knowledge of representations of covering groups has been useful in proving results for reductive groups. Two such instances can be seen in the articles \cite{Shuichiro2014, Shuichiro2015}. These
applications motivate a deeper understanding of genuine representations of
covering groups.  A natural question in this setting is to understand the relation
between a representation and its contragredient, and in particular when the
contragredient can be realized on the same underlying space of the representation.
It is precisely this question that we answer in the present work, for the $n$-fold covers of $\mathrm{GL}(r,F)$.\\

The paper is organized as follows. In Section~\ref{preliminaries}, we briefly recall the Kazhdan-Patterson covers of $G$ and record a few facts that we need. We recall a result of Harish-Chandra about the character of an admissible representation and give references for analogous results in the setting of covering groups. We prove some preliminary results to describe the centralizer of regular semisimple elements in $G$ which is essential in the proof of the main results. We mention a few results of Kable about the lifts of an automorphism. In Section~\ref{Lift of the standard involution}, we discuss some results related to the lifts of the standard involution to $(\widetilde{G}, c)$. In Section~\ref{main results}, we prove the three main results for the cover $(\widetilde{G},c)$. In  Section~\ref{Twisted covers}, we adapt the methods of Section~\ref{main results} to the twisted Kazhdan-Patterson covers $(\widetilde{G}^{(m)}, c^{(m)})$ and study the dualizing property of the lifts of the standard involution.

\section{Preliminaries}\label{preliminaries}

In this section, we set up some preliminaries which we need and recall a few results that will be used throughout this paper.

\subsection{Kazhdan-Patterson covers of $\GL(r)$}\label{Kazhdan-Patterson covers of GL(r)}

Let $F$ be a non-Archimedean local field of characteristic zero. We write $\mathfrak{o}$ for the ring of integers in $F$, $\mathfrak{p}$ for the unique maximal ideal in $\mathfrak{o}$, $\varpi$ for the generator of $\mathfrak{p}$ and $k_{F}$ for the finite residue field of characteristic $p$. For a positive integer $n$, let $\mu_{n}$ be the group of $n$-th roots of unity in $F^{\times}$ and suppose that $|\mu_{n}|=n.$ We let \[\langle \, , \, \rangle_F: F^{\times}\times F^{\times}\rightarrow \mu_{n}\] denote the $n$-th order Hilbert symbol. It 
is known that $\langle \, , \, \rangle_{F}$ is bimultiplicative and defines a non-degenerate pairing on $F^\times/(F^\times)^{n}$ and satisfies
$$\langle a,b\rangle_F \langle b,a\rangle_F= \langle a,1-a\rangle_F=1.$$ 

Let $G=\GL(r,F)$ and $\alpha\in\coh^{2}(G, \mu_{n})$ be a 2-cocycle on $G$ with values in $\mu_{n}$. Associated to the $2$-cocycle $\alpha$,
we can construct a central extension $\widetilde{G}_\alpha$ of $G$ by $\mu_n$. To be precise, take $$\widetilde{G}_\alpha = G \times \mu_n = \{(g, \epsilon): g \in G, ~ \epsilon \in \mu_n\}$$ and define a multiplication on $\widetilde{G}_\alpha$ via the $2$-cocycle $\alpha$ as follows:
$$(g, \epsilon) (h, \eta) = (gh, \epsilon \eta  \alpha(g,h)).$$
It is easy to see that the group $\widetilde{G}_\alpha$ constructed above, fits in the exact sequence given by:
$$1 \longrightarrow \mu_n \xlongrightarrow{i} \widetilde{G}_\alpha \xlongrightarrow{pr} G \longrightarrow 1.$$
Here the embedding $i$ and the projection $pr$ are given by $$i(\epsilon) = (1, \epsilon) \quad \text{and} \quad pr(g, \epsilon) = g.$$
We also have $i(\mu_{n})\subset Z(\widetilde{G}_\alpha)$, where $Z(\widetilde{G}_{\alpha})$ is the center of $\widetilde{G}_{\alpha}$. \\
 
In \cite{KazPat[2]}, Kazhdan and Patterson construct a central extension $\widetilde{G}^{(0)}$ of $G$ with respect to a certain $2$-cocycle $c^{(0)} \in \coh^{2}(G, \mu_{n})$.
For $m \in \mathbb{Z}/n\mathbb{Z}$, they define $2$-cocycles 
$c^{(m)} \in \coh^{2}(G, \mu_{n})$ via,  
\begin{equation*}\label{sigma and sigma c are related}
c^{(m)}(g, h) = c^{(0)}(g, h)  \langle \det(g), \det(h) \rangle_F^m, \quad g, h \in G
\end{equation*}
and describe other ``twisted" extensions $\widetilde{G}^{(m)}$ of $G$. By the Kazhdan-Patterson $n$-fold cover of $G$, we mean the central extension $\widetilde{G}^{(m)}$ of $G$ by $\mu_{n}$ determined by the $2$-cocycle $c^{(m)}$. When $m=0$, for simplicity, we write $c^{(0)}=c$ and $\widetilde{G}^{(0)}$ as $\widetilde{G}$. In what follows, we fix the cover $\widetilde{G}$ with respect to the $2$-cocyle $c$. \\

For $\lambda \in F^{\times}$, write $z_{\lambda} = \lambda 1_{G}.$ Let $Z = \{z_{\lambda}: \lambda^{r-1} \in (F^{\times})^{n} \}.$ Then $Z(\widetilde{G})=pr^{-1}(Z)$ is the center of $\widetilde{G}$ (see Proposition 0.I.1 in \cite{KazPat[2]}). Since $G$ and $\mu_{n}$ are locally compact topological groups, Mackey's theorem (see Theorem 2 in \cite{Mac}) implies that $\widetilde{G}$ is a locally compact topological group and defines a topological central extension of $G$ by $\mu_{n}$. 
It can be shown that the topology on $\widetilde{G}$ has a neighborhood base at the identity consisting of compact open subgroups (see Lemma 3 in \cite{Kab}). Before we give the construction of this basis, we recall a few preliminaries. \\ 

A map $\ell: G\rightarrow \widetilde{G}$ is called a \textit{section} if $pr\circ\, \ell = 1_{G}$. Given a subgroup $S$ of $G$, we say that $\widetilde{G}$ \textit{splits} over $S$, if there exists a homomorphism $s: S\rightarrow \widetilde{G}$ such that $pr\circ\, s =1_{S}$. Let $\ell: G\rightarrow \widetilde{G}$ be the map $\ell(g)=(g,1)$. Then $\ell$ is a section and is called the \textit{natural} or \textit{preferred} section. \\

Let $K_{0}=\GL(r, \mathfrak{o})$ be the maximal compact subgroup in $G.$ For $\lambda\geq 1$, let $K_{\lambda}= 1+\mathfrak{p}^{\lambda}\M(r,\mathfrak{o})$. We know that $\{K_{\lambda}\}_{\lambda\geq 0}$ is a neighborhood base at the identity element in $G$ consisting of compact open subgroups. We can use this base to define a neighborhood base at the identity in $\widetilde{G}$. For a suitable choice of $\lambda_{0}$, it can be shown that $\widetilde{G}$ splits over $K_{\lambda}$, where $\lambda\geq \lambda_{0}$ (see Proposition 0.I.2 in \cite{KazPat[2]}). If $\gcd(n, p)=1,$ we can take $\lambda_0 = 0,$ i.e., $K_{\lambda_{0}} = K_0.$ Let $\kappa: K_{\lambda_{0}}\rightarrow \mu_n$ be a map such that $\kappa^*: K_{\lambda_{0}}\rightarrow \widetilde{G}$ given by $x \mapsto (x,\kappa(x))$ is a splitting of $K_{\lambda_{0}}.$ For $\lambda \geq \lambda_{0}$, let $K^{*}_{\lambda}= \kappa^*(K_\lambda)$.  The collection of compact open subgroups $\{K^{*}_{\lambda}\}_{\lambda\geq \lambda_0}$ forms a neighborhood base at the identity in $\widetilde{G}$.

\begin{remark}\label{twisted covers are l-groups} It can be shown that for $m\in \mathbb{Z}/n\mathbb{Z}$, the other Kazhdan-Patterson covers $(\widetilde{G}^{(m)}, c^{(m)})$ are also topological central extensions of $G$ by $\mu_{n}$ and their topology is given by a neighborhood base at the identity element in $\widetilde{G}$ consisting of compact open subgroups. We refer to Proposition 0.I.2 in \cite{KazPat[2]}. 
\end{remark}

\subsection{Character of an admissible representation}

Let $F$ be a non-Archimedean local field of characteristic zero and $G=G(F)$ be a connected reductive algebraic group defined over $F$. We let $(\pi,V)$ be an irreducible smooth complex representation of $G$. It is known that such representations are admissible. For an admissible representation $(\pi,V)$, we can define a suitable notion of a character. Before we proceed further, we set up some notation and recall a result of Harish-Chandra. \\

Throughout, we let $G=G(F)$ and $(\pi,V)$ to be an irreducible smooth representation of $G$. We let $C^{\infty}_{c}(G)$ denote the space of locally constant complex valued functions on $G$ with compact support. For $f\in C^{\infty}_{c}(G)$, we let $\pi(f): V\rightarrow V$ be the linear operator given by
\[\pi(f)v= \int_{G}f(g)\pi(g)vdg, \quad v\in V,\]
where the integral is with respect to a Haar measure $\mu_{G}$ on $G$ which we fix. If $(\pi,V)$ is an admissible representation, then the trace of the operator $\pi(f)$ is finite for all $f\in C^{\infty}_{c}(G)$. The resulting linear functional
\[\Theta_{\pi}: C^{\infty}_{c}(G)\longrightarrow \mathbb{C}\] given by
\[\Theta_{\pi}(f)=\tr(\pi(f))\]
is called the \textit{distribution character} of $\pi$. It determines the irreducible representation $\pi$ up to equivalence, i.e., 
if $\Theta_{\pi_{1}}(f)=\Theta_{\pi_{2}}(f)$, $\forall f\in C^{\infty}_{c}(G)$, then $\pi_{1}\simeq \pi_{2}$. \\

Let $G_{reg}$ be the subset of regular semisimple elements in $G$. It is known that $G_{reg}$ is an open and dense subset of $G$ whose complement has measure zero. We now state a theorem of Harish-Chandra which is used to define the character of the representation $\pi$.\\ 

\begin{theorem}[Harish-Chandra]\label{regularity theorem} There exists a locally integrable complex valued function $\Theta_{\pi}$ on $G$ such that $\Theta_{\pi}|_{G_{reg}}$ is a locally constant function on $G_{reg}$ and satisfies
\[\Theta_{\pi}(f)
 =\int_{G}f(g)\Theta_{\pi}(g)dg, \quad \forall f\in C^{\infty}_{c}(G).\]
\vspace{0.2 cm}\\
Also, for $x\in G_{reg}, y\in G$, we have $$\Theta_{\pi}(yxy^{-1})=\Theta_{\pi}(x).$$
\end{theorem}
\vspace{0.2 cm}

We refer to \cite{Har} for a proof of the above result. The locally constant function $\Theta_{\pi}$ on $G_{reg}$ is called the \textit{character} of $\pi$. \\

Let $\widetilde{G}$ be a locally compact topological central extension of $G$ by $\mu_{n}$. Let $\xi: \mu_{n}\rightarrow \mathbb{C}^{\times}$ be any injective character of $\mu_{n}$. A representation $(\pi,V)$ of $\widetilde{G}$ is called a \textit{genuine} representation, if for $\epsilon\in \mu_{n}$, $g\in \widetilde{G}$, we have
\[\pi(\epsilon g)=\xi(\epsilon)\pi(g).\]

\vspace{0.2cm}
 The above result of Harish-Chandra also holds in this setting. To be precise, the distribution character of an irreducible admissible genuine representation $\pi$ of $\widetilde{G}$ is represented by a locally integrable function $\Theta_{\pi}$ on $\widetilde{G}$ which is locally constant on $\widetilde{G}_{reg}= pr^{-1}(G_{reg})$ (here $pr: \widetilde{G}\rightarrow G$ is the projection) and satisfies
\[\Theta_{\pi}(y^{-1}xy)=\Theta_{\pi}(x), \textrm{ for } x\in \widetilde{G}_{reg}, y\in \widetilde{G}.\]

We call $\widetilde{G}_{reg}$ as the subset of regular semisimple elements in $\widetilde{G}$. We refer to Theorem 4.3.2 and Corollary 4.3.3 in \cite{Wen} where results about the character theory for metaplectic groups are discussed. See also Theorem I.5.1 in \cite{KazPat[2]} where the result is only mentioned.

\subsection{Centralizer of regular semisimple elements in $\GL(r,F)$}\label{centralizer discussion}

Let $F$ be a non-Archimedean local field of characteristic zero and 
let $G=\GL(r,F)$. Given a finite dimensional algebra $K/F$, we let $N_{K/F}$ denote the norm of $K/F$. In this section, we give a description of the centralizer of a regular semisimple element in $G$. Realizing the centralizer of a regular semisimple element as a certain subset of the algebra $\M(r,F)$, we also discuss some of its properties under the norm map.  \\

Recall that a matrix $g\in G$ is said to be regular semisimple, if its characteristic polynomial $p_g(X)$ has no repeated roots
over $\overline F$.  Equivalently, the minimal polynomial $m_g(X)$ is a product of distinct monic irreducible polynomials over $F$ and has degree $r$. In this case, one has $m_g(X)= p_g(X).$ Let $g\in G_{reg}$. Consider the subalgebra $F[g]$ of $\M(r,F)$.   
The evaluation map $$f(X) \longmapsto f(g): F[X]\longrightarrow F[g]$$ has kernel
$(m_g(X)),$ and it follows that
\[
F[X]/(m_g(X)) \simeq F[g].
\]
In particular, we have $\dim_F F[g]=r.$ Suppose that   
\[
p_g(X)= \prod_{i=1}^{s} p_i(X),
\]
where $p_{i}(X)$'s are the distinct monic irreducible polynomials over $F$. Let
\begin{equation*}
F_i:=F[X]/(p_i(X)) \quad \text{and} \quad d_i:=\deg (p_i(X)) =[F_i:F].
\end{equation*}

\begin{lemma}\label{lem:centralizer-Fg} For $g \in G_{reg}$, we have
\[
F[g] = \{A\in \M(r,F):Ag=gA\}.
\]
Consequently,
\[
Z_{G}(g)=F[g]^{\times}.
\]
\end{lemma}

\begin{proof} Since $g$ commutes with every element in $F[g]$, we have $F[g]^{\times}\subset Z_{G}(g)$. Let $A\in \M(r,F)$ such that $Ag=gA$. Since $m_g(X)=p_g(X)$, the matrix $g$ is cyclic, i.e., there exists $v\in F^r$
such that
\[
\{v,gv,\ldots,g^{r-1}v \}
\]
is an $F$-basis of $F^r$.
Since
\[
\{v,gv,\ldots,g^{r-1}v \}
\]
is a basis, there is a unique polynomial $f(X) \in F[X]$ of degree less than $r$ such
that
\[
Av=f(g)v.
\]
The relation $Ag=gA$ implies that $Ag^{k}=g^{k}A$ for every $k\geq 0$.  Hence, for
$0\leq k<r$,
\begin{align*}
A(g^{k}v)
&=g^{k}(Av)=g^{k}f(g)v=f(g)(g^{k}v).
\end{align*}
Thus $A$ and $f(g)$ agree on the basis
$\{v,gv,\ldots,g^{r-1}v \}$, and therefore $A=f(g)$.  Passing to invertible
elements, gives $Z_G(g)\subset F[g]^{\times}$. Hence the result.
\end{proof}

\begin{lemma}\label{identifying F[g]^x as a product of Fi's} For $g \in G_{reg}$, we have   
$$F[g]\simeq F_1\times\cdots\times F_s.$$ 
In particular, $d_{1}+\cdots+d_{s}=r.$ Consequently,
\[
Z_{G}(g)=F[g]^{\times}\simeq F_{1}^{\times} \times \cdots \times F_{s}^{\times}.
\]
\end{lemma}

\begin{proof} Since the ideals $(p_i(X))$ are pairwise comaximal, the Chinese Remainder
Theorem gives
\[
F[g]\simeq F_{1} \times \cdots \times F_{s}. 
\]
Passing to the units, we get
\[
Z_{G}(g)=F[g]^{\times}\simeq F_{1}^{\times} \times \cdots \times F_{s}^{\times}. \qedhere
\]
\end{proof}

\begin{remark} Let  $u = (u_{1},\ldots ,u_{s}) \in F[g]^\times.$ Clearly, we have $u\in (F[g]^{\times})^{2}$ if and only if each
$u_{i} \in (F_{i}^{\times})^{2}$.
\end{remark}

\begin{lemma}\label{lem:det-equals-norm} For $u=(u_{1}, \ldots, u_{s})\in F[g]^{\times}$, let \begin{equation*}\label{definition of N(u)}N(u):=\prod_{i=1}^{s} N_{F_{i}/F}(u_{i}).\end{equation*}
Under the embedding $F[g]\hookrightarrow \M(r,F)$, we have
\[
\det(u)=N(u).
\]
\end{lemma}

\begin{proof} Choose a cyclic vector $v\in F^{r}$. The map
\[
\Phi:F[g]\longrightarrow F^{r}; \quad a\longmapsto av,
\]
is an $F$-linear isomorphism. Let $m_{u}\in \End_{F}(F[g])$ denote the map given by multiplication by $u$, and $T_u\in \End_{F}(F^{r})$ denote the map given by the usual action of $u$ on $F^{r}.$ For $a\in F[g]$, we have
\begin{align*}
\Phi(m_u(a))
&=\Phi(ua)=(ua)v=u(av)=T_u(\Phi(a)).
\end{align*}
Therefore, 
\begin{equation*}
T_u=\Phi\circ m_u\circ\Phi^{-1}, \textrm{ and } \det(u) = \det(m_u).
\end{equation*}
Under the identification $\displaystyle F[g]\simeq  \prod\limits_{i=1}^{s} F_{i} $, multiplication by $u$ is the direct sum of the maps
\[
m_{u_i}:F_{i}\longrightarrow F_{i}; \quad a_{i}\longmapsto u_{i}a_{i}.
\]
Thus,
\begin{align*}
\det(u) &=  \det(m_u)=\prod_{i=1}^s {\det}(m_{u_i})=\prod_{i=1}^s N_{F_i/F}(u_i)=N(u). \qedhere
\end{align*}
\end{proof}

Let $n\geq 2$ be an integer. Denote by 
$\langle \, \  , \,  \rangle_{K}$ and $\langle \, \  , \,  \rangle_{F}$ the $n$-th order Hilbert symbol on $K$ and $F$ respectively. We record an important property relating $\langle \, \  , \,  \rangle_{K}$ and $\langle \, \  , \,  \rangle_{F}$ in the following lemma. 

\begin{lemma}[{\cite[p.~48]{KazPat[2]}}]
\label{relation between Hilbert symbol on K and F}
For every $a\in F^\times$ and $b\in K^\times$, we have
\begin{equation*}
\langle a,b\rangle_{K}
    =
    \left\langle a,N_{K/F}(b)\right\rangle_{F}.
\end{equation*}
\end{lemma}

\subsection{Some results about lifts of an automorphism of $G$ to $\widetilde{G}$}\label{results we need}

We recall a few results from \cite{Kab} which we use in our main result. Let $\widetilde{G}$ be a central extension of a group $G$ by an abelian group $A$. Let $pr: \widetilde{G}\rightarrow G$ be the projection map, $s: G\rightarrow \widetilde{G}$ be a section of $pr$ and $\alpha$ be the 2-cocycle representing the class of this central extension in $\coh^{2}(G,A)$ with respect to the section $s$. If $f: G\rightarrow G$ is an automorphism  of $G$, then a lift of $f$ is an automorphism  $\tilde{f}: \widetilde{G}\rightarrow \widetilde{G}$ such that 
\begin{equation}\label{lift-condition}
    pr(\tilde{f}(g))=f(pr(g)), \forall g\in \widetilde{G} \quad \text{and} \quad \tilde{f}(a) = a, \forall a\in A.
\end{equation}

\begin{remark}
  The corresponding statement for anti-automorphisms $f$ and $\tilde{f}$ is obtained by replacing the second condition in \eqref{lift-condition} with  $$ \tilde{f}(a) = a^{-1}, \quad \forall a \in A. $$ 
\end{remark} 
The group $\Aut(G)$ acts on $\coh ^{2}(G,A)$ by $$ \phi [\Lambda]= [\Lambda \circ (\phi^{-1}\times \phi^{-1})].$$

\begin{proposition}\label{description of the set of liftings} Let $\mathcal{L}(f)$ denote the set of all lifts of $f$. We have the following: 
\begin{enumerate}
\item[\emph{1)}] The set $\mathcal{L}(f)$ is non-empty if and only if $f [\alpha]=[\alpha]$.
\smallskip
\item[\emph{2)}] If $\mathcal{L}(f)$ is non-empty, then $\mathcal{L}(f)$ is a principal homogeneous space for the group $\Hom(G,A)$ under the action
\[(\phi.\tilde{f})(g)= \phi(pr(g))\tilde{f}(g).\]
\end{enumerate}
\end{proposition}

We also need the following result (see Corollary 1 in \cite{Kab} for a proof) where the continuity properties of the lift in the case when $G=\GL(r,F)$ is discussed. 

\begin{proposition}\label{continuity of lift} Let $F$ be a non-Archimedean local field and suppose that the group of $n$-{th} roots of unity in $F$ has order $n$. Let $\langle \, , \, \rangle_F$ be the $n$-th order Hilbert symbol on $F$ and $\widetilde{G}$ the corresponding metaplectic group. Then the lift of any topological automorphism (resp. anti-automorphism) of $G$ to $\widetilde{G}$ is also a topological automorphism (resp. anti-automorphism).
\end{proposition}

\section{Lifts of the standard involution of $G$ to $\widetilde{G}$}\label{Lift of the standard involution}

We continue with the same notations as in the earlier sections. Fix a positive integer $r \geq2$. Let $G=\GL(r,F)$ and $S=\SL(r,F)$. Let $\widetilde{G}$ be the $n$-fold metaplectic cover of $G$ with respect to the $2$-cocycle $c\in \coh^{2}(G, \mu_{n})$. Let $H$ be the subgroup of diagonal matrices in $G$. For $g\in G$, we write $^{\top}g$ for the transpose of the matrix $g$.  Let $\theta$ be the involution on $G$ defined by $g \mapsto ~^{\top}g$. Let $\tau$ be the standard involution on $G$ given by $g \mapsto w_{0}{}^\top gw_{0},$ where $w_{0}\in G$ is the matrix with anti-diagonal entries one. In this section, we discuss some results about the lifts of $\tau$ to $\widetilde{G}$, which we will need in proving the main results. We begin by defining a certain lift $\sigma$ of $\tau$ and show that it is an involutive anti-automorphism of $\widetilde{G}$.
Throughout, we let $\beta\in \Aut(G)$ to be the automorphism given by $g \mapsto ~^{\top}g^{-1}.$ 

\begin{lemma}\label{c beta and c are cohomologous on G} Let $c^\beta: G \times G \to \mu_n$ be the map
$$c^\beta(x, y) =c(\beta(x), \beta(y)).$$ 
Then $c^\beta$ and $c$ represent the same cohomology class in $\coh^2(G, \mu_n).$
\end{lemma}
We refer to  \cite[Proposition~3.1]{KazPat} and \cite[Lemma~2 \& Proposition~1]{Kab} for a proof. \\

By Lemma~\ref{c beta and c are cohomologous on G}, we can choose a map $\delta:G \to \mu_n$ such that  
\begin{equation}\label{sigma beta and sigma cohologous main equation}
c^{\beta}(x, y) =c(x, y) \frac{\delta(xy)}{\delta(x) \delta(y)}, \forall x, y \in G.
\end{equation}
In fact, we can choose $\delta$ satisfying \eqref{sigma beta and sigma cohologous main equation} such that $\delta|_{H} =1,$ in which case it is unique. In what follows, we fix such a choice of $\delta.$

\begin{lemma}\label{to show beta tilde is an involution}
    Let $\phi:G \to \mu_n$ be the map defined as follows:
    $$\phi(x) = \delta(\beta(x)) \delta(x).$$ Then
    \begin{enumerate}[label=\emph{\arabic*)}]
        \item $\phi$ is a character of $G.$
\vspace{0.1 cm}
        \item $\phi =1.$
    \end{enumerate}
\end{lemma}

\begin{proof}
Using $\beta^{2} =1$ together with \eqref{sigma beta and sigma cohologous main equation}, the proof of 1) follows. Since the commutator subgroup of $G$ is $S,$ every character of $G$ factors through the determinant map. Hence, any character of $G$ with values in $\mu_n$ is of the type $x \mapsto \chi(\det(x)),$ where $\chi:F^\times \to \mu_n$ is a character. We know that every character $\chi$ of $F^\times$ has the form $\chi(a) = \langle \gamma, a \rangle_F$ for some $\gamma \in F^\times.$ Thus, $$\phi(x) = \langle \gamma, \det(x) \rangle_F.$$ For $a \in F^\times,$ consider $x = diag(a,1, \dots, 1) \in H.$ Since $\delta|_{H} =1$ and $\beta$ fixes $H,$ we have  $\phi(x) = 1.$ Consequently, 
$$\langle \gamma, a \rangle_F =1 \quad \text{for every} \quad a \in F^\times.$$ 
The non-degeneracy of the Hilbert symbol implies that  $ \gamma \in (F^{\times})^n.$ The result follows.
\end{proof}

\begin{lemma}\label{lift of the automorphism beta} Let $\tilde{\beta}: \widetilde{G} \to \widetilde{G}$ be the map given by 
\begin{equation*}\label{definition of tilde beta}   
\tilde{\beta} (x, \epsilon) = (\beta(x), \epsilon \delta(x)).\end{equation*} Then
\begin{enumerate}[label=\emph{\arabic*)}]
    \item $\tilde{\beta}$ is an automorphism of $\widetilde{G}.$
\vspace{0.1 cm}
    \item $\tilde{\beta}^2 =1.$ 
\vspace{0.1 cm}
    \item $\tilde{\beta}$ is a lift of $\beta.$
\end{enumerate}
\end{lemma}

\begin{proof} Using \eqref{sigma beta and sigma cohologous main equation}, part 1) follows immediately. By Lemma~\ref{to show beta tilde is an involution}, part 2) follows. Indeed,
$$\tilde{\beta}(\tilde{\beta}(x, \epsilon)) = (\beta^2(x), \epsilon \delta(\beta(x)) \delta(x)) = (x, \epsilon).$$
For 3), it is easy to verify that $\beta \circ pr = pr \circ\, \tilde{\beta}.$
\end{proof}

\begin{lemma}\label{definition of the involution rho} Let $\rho: \widetilde{G} \to \widetilde{G}$ be the map given by \[\rho(g) \mapsto \tilde{\beta}(g^{-1}), \forall g=(x, \epsilon)\in \widetilde{G}.\] Then
\begin{enumerate}[label=\emph{\arabic*)}]
 \item $\rho(x , \epsilon) = (\theta(x), \epsilon^{-1} \delta(x^{-1}) c(x, x^{-1})^{-1}).$
 \vspace{0.1 cm}
    \item $\rho$ is an anti-automorphism of $\widetilde{G}$ and $\rho^{2}=1.$ 
\vspace{0.1 cm}
    \item $\rho$ is a lift of  $\theta$.
\end{enumerate}
\end{lemma}
The proof of the above lemma follows from Lemma~\ref{lift of the automorphism beta}.

\begin{lemma}\label{definition of the involution sigma} Let $\sigma: \widetilde{G} \to \widetilde{G}$ be the map given by 
\[\sigma(g)=(w_{0},1) \rho(g) (w_{0},1)^{-1}, \forall g\in \widetilde{G}.\]
Then
\begin{enumerate}[label=\emph{\arabic*)}]
\item $\sigma$ is an anti-automorphism of $\widetilde{G}.$
\vspace{0.1 cm}
\item $\sigma^2 =1.$ 
\vspace{0.1 cm}
\item $\sigma$ is a lift of $\tau.$
\end{enumerate}
\end{lemma}

\begin{proof} The proof of 1) and 3) follows from straightforward computations using Lemma~\ref{definition of the involution rho}. For 2), let $g \in \widetilde{G}.$ We have, 
\begin{equation*}\label{sigma is an involution 1}
\sigma(\sigma(g)) = (w_0,1) \rho((w_0, 1)^{-1}) g \rho((w_0, 1)) (w_0, 1)^{-1}.
\end{equation*}
Note that $w_{0}^{2} = 1_{G}$ and $\beta(w_{0}) = w_{0}$. Thus, applying part 1) of Lemma~\ref{definition of the involution rho}, we get
$$\sigma(\sigma(g)) = g (1_{G}, \delta(w_{0})^{2}).$$
It follows from \eqref{sigma beta and sigma cohologous main equation}  that $\delta(w_{0})^{2} =1.$ Hence the result. 
\end{proof}

By Proposition 1 in \cite{Kab}, there are several lifts of the standard involution $\tau$ to $\widetilde{G}$. In the following lemma, we discuss the other lifts of $\tau$ to $\widetilde{G}$.

\begin{lemma}\label{other lifts of tau}  For $g = (x, \epsilon) \in \widetilde{G},$ we write $\triangle(g) = \det(x).$ For  $\alpha \in F^{\times},$ let $\sigma_{\alpha}:  \widetilde{G} \to \widetilde{G}$ be the map given by
\begin{equation*}\label{definiton of sigma-alpha}
\sigma_{\alpha}(g)= \langle \alpha, \triangle(g)\rangle_{F} \,\sigma(g).\end{equation*} 
Then 
\begin{enumerate}[label=\emph{\arabic*)}]
\item $\sigma_{\alpha}$ is an anti-automorphism on $\widetilde{G}.$
\vspace{0.1 cm}
\item $\sigma_{\alpha}^2 =1.$ 
\vspace{0.1 cm}
\item $\sigma_{\alpha}$ is a lift of $\tau$. 
\end{enumerate}
\end{lemma}

\begin{proof} The proof follows from Lemma~\ref{definition of the involution sigma}. Parts 1) and 3) are immediate by direct computations. We prove only 2).
Write $h = \sigma_{\alpha}(g).$ Then
\begin{align*}
\sigma_{\alpha}(\sigma_{\alpha}(g))=\sigma_{\alpha}(h) & = \langle \alpha, \triangle(h)\rangle_{F} \,  \sigma(h) \\
    & = \langle \alpha, \triangle(h)\rangle_{F} \langle \alpha, \triangle(g)\rangle_{F}^{-1}\, \sigma^{2}(g)
\end{align*}
Since $\triangle(h) =  \triangle(g)$ and $\sigma^{2} =1,$ the result follows.  \end{proof}

\begin{remark} Note that any character $\phi:G \to \mu_n$ is given by $x \mapsto \langle \alpha, \det(x) \rangle_F$ for some $\alpha\in F^{\times}$. From Proposition~\ref{description of the set of liftings}, we have that any lift of $\tau$ to $\widetilde{G}$ is of the form $\sigma_{\alpha}$ for $\alpha\in F^{\times}$. Moreover, $\sigma_{\alpha} = \sigma_{\gamma}$ if and only if $\alpha \gamma^{-1} \in (F^{\times})^n$. Thus, $\{ \sigma_{\alpha}: \alpha \in F^{\times}/(F^{\times})^n \}$ is the set of all distinct lifts of $\tau$ to $\widetilde{G}.$ 
\end{remark}

\section{Proof of the main results}\label{main results} In this section, we prove the three main results of this paper. Throughout, we let $F$ be a non-Archimedean local field of characteristic zero and $k_{F}$ denote the finite residue field of characteristic $p$. Let $r\geq 2$ be a positive integer and let $G=\GL(r,F)$. Let $\theta(g)={}^{\top}g$ and $\tau(g)=w_{0}{}^{\top}g w_{0}$. We let $\widetilde{G}$ be the $n$-fold Kazhdan-Patterson cover of $G$ with respect to the $2$-cocycle $c$. Let $\rho$ be a lift of $\theta$ to $\widetilde{G}$. We write $\mu_{\widetilde{G}}$ for the Haar measure on $\widetilde{G}$ and $\Aut_{c}(\widetilde{G})$ for the group of continuous automorphisms of $\widetilde{G}$.  $\Rep(\widetilde{G})$ denotes the category of representations of $\widetilde{G}$.
\\

It is known that for $r=2$ and $n\geq 3$, any lift of the standard involution $\tau$ to $\widetilde{G}$ is not a duality involution (see \cite{Kumar-Renu-Sanjeev-Varsha}). This motivates the first question considered in this article. In our first main result, assuming that $r\geq 2$ and $n\geq 3$,  we show that any lift of $\tau$ to the $n$-fold cover $\widetilde{G}$ is never a duality involution. This negative result motivates the other questions, where we focus exclusively on understanding the dualizing property of the lifts only for the $2$-fold cover of $G$. In our second main result, under the assumption that $r=2t$ is an even integer, we show that any lift of $\tau$ to the $2$-fold cover $\widetilde{G}$ is always a duality involution. When $r=2t+1$ is an odd integer, as our third main result, we show that there is a unique dualizing lift of $\tau$ to the $2$-fold cover $\widetilde{G}$. Before we continue, we record a few results that we need.

\begin{lemma}\label{sigma-rho dualizing property} Let $\rho$ be any lift of $\theta$ to $\widetilde{G}$ and let $\sigma$ be defined as in Lemma~\ref{definition of the involution sigma}. We have $\rho$ is dualizing if and only if $\sigma$ is dualizing. 
\end{lemma}
 
\begin{proof}  Let $(\pi,V)\in \Rep(\widetilde{G})$ such that it is irreducible and admissible. Suppose that $\pi^{\rho}\simeq \pi^{\vee}$. Consider the map $L: V \to V^\vee$ given by $L(v)=T(\pi(\tilde{w}_0)v)$, where $ \tilde{w}_0 = (w_0, 1)$ and $0 \neq T\in \Hom_{\widetilde{G}}(\pi^\rho, \pi^\vee)$. It is easy to see that $L\in \Hom_{\widetilde{G}}(\pi^\sigma, \pi^\vee)$ and is an isomorphism. The converse is similar.  
\end{proof}

\begin{lemma}\label{results-on-haarmeasure} Let $\gamma\in \Aut_{c}(\widetilde{G})$. Then
\begin{enumerate}
\item[\emph{1)}] There exists $c_{\gamma}>0$ such that $\mu_{\widetilde{G}}\circ \gamma = c_{\gamma}\mu_{\widetilde{G}}.$ 
\smallskip
\item[\emph{2)}] The map $\gamma\longmapsto c_{\gamma}: \Aut_{c}(\widetilde{G}) \longrightarrow \mathbb{R}^{\times}_{>0}$ is a homomorphism.
\smallskip
\item[\emph{3)}]   If $\gamma^{2}=1,$ then  $\mu_{\widetilde{G}}\circ \gamma =\mu_{\widetilde{G}}.$
\end{enumerate}
\end{lemma}

\begin{proof}\hfill 
\begin{enumerate} 
\item[1)] Write $\nu_{\widetilde{G}}=\mu_{\widetilde{G}}\circ\gamma.$
It is easy to see that $\nu_{\widetilde{G}}$ is a left invariant Haar measure on $\widetilde{G}$. Indeed, for $g\in G$ and $U\subset G$, we have
$$\nu_{\widetilde{G}}(gU)=\mu_{\widetilde{G}}(\gamma(g)\gamma(U))=\mu_{\widetilde{G}}(\gamma(U))=\nu_{\widetilde{G}}(U).$$
Using the uniqueness of Haar measure, the result follows. \\

\item[2)] For an open set $U\subset G$, we have
\begin{align*}
c_{\gamma_1\gamma_2}\mu_{\widetilde{G}}(U)
&= (\mu_{\widetilde{G}}\circ\gamma_1\gamma_2)(U)\\
&= (\mu_{\widetilde{G}}\circ\gamma_1)(\gamma_2(U))\\
&= c_{\gamma_1}(\mu_{\widetilde{G}}\circ\gamma_2)(U)\\
&= c_{\gamma_1}c_{\gamma_2}\mu_{\widetilde{G}}(U).
\end{align*}

\item[3)] Let $K$ be any compact open subset of $\widetilde{G}$. Since $K$ is open, using
properties of the Haar measure it follows that $\mu_{\widetilde{G}}(K)>0.$
Since $\gamma^2=1$, we have
$$c_\gamma\mu_{\widetilde{G}}(K)= \bigl(\mu_{\widetilde{G}}\circ\gamma^{-1}\bigr)(K)= c_{\gamma^{-1}}\mu_{\widetilde{G}}(K)= c_\gamma^{-1}\mu_{\widetilde{G}}(K).$$
It follows that $c_\gamma=1$ and hence the result. Consequently, for any
measurable subset $X\subset\widetilde{G},$  we have
\[
\mu_{\widetilde{G}}\bigl(\gamma(X)\bigr)
=
\mu_{\widetilde{G}}(X).  \qedhere
\]
\end{enumerate}
\end{proof}

\begin{lemma}\label{pi rho f is same as pi f check rho} For $f \in C_{c}^{\infty}(\widetilde{G})$, define $f^{\vee},f^{\rho} \in C_{c}^{\infty}(\widetilde{G})$ as follows:
\[f^{\vee}(g)=f(g^{-1}) \quad \text{and} \quad f^{\rho}(g)=f(\rho(g)).\]
For an irreducible admissible and genuine representation $\pi\in \Rep(\widetilde{G})$, we have the following: 
\begin{enumerate}
\item[\emph{1)}] $\pi^{\rho}(f)=\pi((f^{\vee})^{\rho})$.
\smallskip
\item[\emph{2)}] $\pi^{\vee}(f)= ~^{\top}\pi(f^{\vee})$.
\smallskip
\item[\emph{3)}] $\Theta_{\pi^{\vee}}(g)=\Theta_{\pi}(g^{-1})$, $g\in \widetilde{G}_{reg}$. \\
\end{enumerate}
\end{lemma}

\begin{proof} \hfill 

\begin{enumerate}[label=\arabic*)] 
\item Using a simple change of variables and Lemma~\ref{results-on-haarmeasure} the result follows. Indeed,
 \begin{align*} \pi((f^{\vee})^{\rho}) &= \displaystyle\int_{\widetilde{G}} f(\rho(g^{-1})) \pi(g) dg \\
&= \int_{\widetilde{G}} f(\tilde{\beta}(g)) \pi(g) dg  \\
&= \int_{\widetilde{G}} f(g) \pi(\tilde{\beta}(g)) dg  \\
&=  \pi^{\rho}(f).
\end{align*}

\item For $g \in \widetilde{G},$ let ${}^{\top}{\pi(g)}$ denote the transpose of the operator $\pi(g).$ Since the restriction of the modular character $\delta_{\widetilde{G}}$ to $\mu_n$ is trivial, it follows that $\delta_{\widetilde{G}}$ factors through the projection map $pr: \widetilde{G} \to G$, and we have $\delta_{\widetilde{G}} = \delta_{G} \circ pr.$ Since $G$ is unimodular, we get $\delta_{\widetilde{G}} = 1.$ Thus we have 
\begin{align*}
\pi^\vee(f) &= \displaystyle\int_{\widetilde{G}} f(g) \pi^\vee(g) dg\\
&= \displaystyle\int_{\widetilde{G}} f(g)^\top\pi(g^{-1}) dg\\
&= \displaystyle\int_{\widetilde{G}} f^\vee(g)^\top\pi(g) \delta_{\widetilde{G}}(g^{-1})dg \\ 
& = ~^{\top}\pi(f^\vee). 
\end{align*}

\item Since the trace is invariant under taking transpose, it follows that $$\tr(\pi^\vee(f)) = \tr(\pi(f^\vee)),$$ for all $f \in C_c^\infty(\widetilde G).$ Using Theorem \ref{regularity theorem} for $\widetilde{G}$, the result follows. 
Indeed, for $g\in \widetilde{G}_{reg}$ 
\begin{align*}
\Theta_{\pi^{\vee}}(g) &= \int_{\widetilde{G}_{reg}} f(g)\Theta_{\pi^\vee}(g)dg \\
& = \int_{\widetilde{G}_{reg}} f^\vee(g)\Theta_{\pi}(g)dg \\
& = \int_{\widetilde{G}_{reg}} f(g)\Theta_{\pi}(g^{-1})dg\\
&= \Theta_{\pi}(g^{-1}).
\end{align*}
\end{enumerate}
\end{proof}

\begin{lemma}\label{restriction-mu-n-of-central-character-injective}
	Let $\pi$ be an irreducible admissible genuine representation of $\widetilde{G}.$ Let $\omega_\pi$ be the central character of $\pi.$ Then
	\begin{enumerate}[label=\emph{\arabic*)}]
		\item ${\omega_\pi}\big |_{\mu_n}$ is injective.
		\smallskip
		
		\item $\omega_\pi(\eta^2) = 1,$ for every $\eta\in \mu_n$ if and only if $n=2.$
	\end{enumerate}
\end{lemma}

\begin{proof} Let $\eta\in \mu_n$. Since $\pi$ is a genuine representation,  it follows that $ \omega_\pi(\eta) = \xi(\eta)$ for some injective character $\xi: \mu_n \to \mathbb{C}^\times.$ This proves 1).  If $n=2$, clearly we have $\omega_\pi(\eta^2) = 1.$ Suppose that $\omega_\pi(\eta^2) = 1$, then $\omega_\pi(\eta) \in \{\pm 1\}.$  If $n \geq 3,$ there would exist $\eta_{1}, \eta_{2} \in \mu_{n}$ with $\eta_{1} \neq \eta_{2}$ such that $\omega_{\pi}(\eta_{1}) = \omega_\pi(\eta_{2}),$ contradicting part 1). 
\end{proof}

In what follows, it suffices to show that a lift $\rho$ of $\theta$ to $\widetilde{G}$ is either a non-dualizing or a dualizing involution to prove that the corresponding lift $\sigma$ of $\tau$ is non-dualizing or dualizing.

\subsection{Non-dualizing lifts of $\tau$} \label{Non-dualizing lifts of tau}

Let $F$ be a non-Archimedean local field of characteristic $0$. Let $k_{F}=\mathfrak{o}/\mathfrak{p}$ be the finite residue field of $F$, where $\mathfrak{o}$ is the ring of integers in $F$ and $\mathfrak{p}$ is the unique maximal ideal in $\mathfrak{o}$. Let $G=\GL(r,F)$ and $\widetilde{G}$ be the $n$-fold Kazhdan-Patterson cover of $G$ with respect to $2$-cocycle $c\in \coh^{2}(G, \mu_{n})$ (see Section~\ref{Kazhdan-Patterson covers of GL(r)}; for more details we refer to Chapter~$0$ in \cite{KazPat[2]}. We note that \cite{KazPat[2]} uses the notation $\sigma^{(0)}$ for the $2$-cocycle in place of $c$). In this section, we provide two proofs of the fact that the lifts $\sigma_\alpha$ of the standard involution $\tau$ to $\widetilde{G}$ is not a duality involution. The first proof is constructive (inspired by a comparison of the distribution characters) and works only for the case when $r\geq 3$ and with some restrictions on the residual degree of $F$. To be precise, we assume that $r,n\geq 3$ with $\gcd(n,p)=1$ and construct a compact open subset of $\widetilde{G}$ and a function $f\in C_{c}^{\infty}(\widetilde{G})$ supported on this compact set, such that the distribution characters of $\pi^{\rho}$ and $\pi^{\vee}$  evaluated at $f$ are distinct. The second proof is more elementary (inspired by a comparison of the central characters) and it includes all the cases $r\geq 2$ and $n\geq 3$ and works independently of any restrictions on the residual degree of $F$. 

\subsubsection{Comparison of distribution characters (proof~1):} Assume that $r,n \geq 3$ and $\gcd (n,p)=1$, where $\Char(k_{F})=p$. Let $\beta\in \Aut(G)$ be the involutive automorphism of $G$ defined as $\beta(g)=~^{\top}g^{-1}.$ Let $\delta:G\rightarrow \mu_{n}$ be the unique map such that $\delta|_H=1$, where $H$ is the diagonal subgroup of $G$ and satisfying \eqref{sigma beta and sigma cohologous main equation}: 
\[c^\beta(x,y):= c(\beta(x),\beta(y))= \frac{\delta(xy)}{\delta(x) \delta(y)} \ c(x,y).\] 
Consider the lift $\tilde{\beta}$ of $\beta$ as in Lemma~\ref{lift of the automorphism beta} given by: 
\begin{equation*} \tilde{\beta}(x,\epsilon)= (\beta(x), \epsilon\delta(x)).
\end{equation*}
Let $K_{0}=\GL(r,\mathfrak{o})$. Since $\gcd (p,n)=1$, $\widetilde{G}$ splits over $K_0.$ Let $\kappa:K_0 \to \mu_n$ be a map such that $\kappa^*:K_0 \to \widetilde{G}$ given by $\kappa^*(x)=(x, \kappa(x))$ is a homomorphism and satisfies $pr\circ\, \kappa^*= 1_{K_0}$.
We observe that
$\kappa^*$ is a homomorphism if and only if $\kappa(xy)=\kappa(x)\kappa(y)c(x,y)$ for all $x,y \in K_0.$ Since $\kappa^*$ is a splitting, it follows that $\tilde{\beta} \circ \kappa^*$ and $\kappa^* \circ \beta$ are two lifts of $\beta$ in the following sense: $\tilde{\beta} \circ \kappa^*$ and $\kappa^* \circ \beta$ are maps from $K_0$ to $\widetilde{G}$ satisfying
\[ (pr \circ \tilde{\beta} \circ \kappa^*)(x) = \beta(x)= (pr \circ \kappa^* \circ \beta)(x).\]

\begin{lemma}\label{phi is trivial}
The map $\phi:K_0 \longrightarrow \mu_n$ defined as $$\phi(x):= (\tilde{\beta} \circ \kappa^*)(x)(\kappa^*\circ \beta)(x) ^{-1}$$ is a group homomorphism. In fact, $\phi$ is trivial.
\end{lemma}

\begin{proof} We have
    \begin{align*}
    \phi(xy) = \frac{(\beta(xy), \kappa(xy) \delta(xy))}{(\beta(xy), \kappa(\beta(xy)))}.
    \end{align*}
On the other hand, using $\kappa(xy)=\kappa(x)\kappa(y)c(x,y)$, we get

    \begin{align*}
        \phi(x)\phi(y) 
        &=\frac{(\beta(x), \kappa(x) \delta(x))}{(\beta(x), \kappa(\beta(x)))} \frac{(\beta(y), \kappa(y) \delta(y))}{(\beta(y), \kappa(\beta(y)))}\\[2pt]
        &= \frac{(\beta(xy), \kappa(x)\kappa(y) \delta(x) \delta(y) c^\beta(x,y))}{(\beta(xy), \kappa(\beta(x))\kappa(\beta(y)) c^\beta(x,y))}\\[2pt]
        &= \frac{(\beta(xy), c(x,y)^{-1} \kappa(xy) \delta(x) \delta(y) c^\beta(x,y))}{(\beta(xy), \kappa(\beta(xy)) c^\beta(x,y)^{-1}c^\beta(x,y))}\\[2pt]
        &=\frac{(\beta(xy), \kappa(xy) \delta(xy))}{(\beta(xy), \kappa(\beta(xy))}.
    \end{align*}
Hence $\phi$ is a homomorphism. Further,  as $r \geq 3,$  the commutator $[K_0,K_0]=\SL(r, \mathfrak{o})$. 
Since $\gcd (p,n)=1$, we have $\langle \ , \ \rangle_{F} ~\big |_{\mathfrak{o}^{\times} \times \mathfrak{o}^{\times}}=1$. Thus, by Proposition 3.1 in \cite{KazPat}, it follows that $\phi$ is trivial.
\end{proof}
For $m \geq 1$, let $K_m=1+\mathfrak{p}^{m}\M(r, \mathfrak{o})$. Clearly, $K_{m}$ is a subgroup of $K_0$. Since $\phi$ is trivial on $K_{0}$, we have $\phi|_{K_{m}}$ is trivial for all $m \geq 1$. 

\begin{lemma}\label{betatildefixesKm} For $m \geq 0,$ write $K_m^* = \kappa^*(K_m)$. We have $$\tilde{\beta}(K_m^*)=K_m^*.$$
\end{lemma}

\begin{proof} Using Lemma~\ref{phi is trivial} and $\beta(K_{m})=K_{m}$, the proof follows. Indeed,
   \begin{align*}
       \tilde{\beta}(K_m^*) &= \tilde{\beta}(\kappa^*(K_m))\\
       &= \kappa^*(\beta(K_m))\\
       &= \kappa^*(K_m)\\
       &= K_m^*. 
   \end{align*}
\end{proof}

For $A\subset \widetilde{G}$, we let $\chi_{A}$ denote the characteristic function of $A$. Let $\rho:\widetilde{G}\rightarrow \widetilde{G}$ be the involutive anti-automorphism on $\widetilde{G}$ defined in Lemma~\ref{definition of the involution rho}:  
\[\rho(g)=\tilde{\beta}(g^{-1}), \forall g \in \widetilde{G}.\]
Throughout, we write $A^{\rho}=\rho(A)$ for the image under $\rho$ of 
a subset $A$ of $\widetilde{G}$. We write
 $\mu=\mu_{\widetilde{G}}$ for the Haar measure on $\widetilde{G}$. For $K$ a compact open subgroup of $\widetilde{G}$, we let $e_K\in C_{c}^{\infty}(\widetilde{G})$ be the function 
\[e_K(g) = \begin{cases}
	\mu(K)^{-1}, & g \in K \\
	0, & g \notin K.
\end{cases}\]

\begin{theorem}[Non-Dualizing lifts]\label{Non-Dualizing lift} Let $(\pi, V)$ be an irreducible admissible genuine representation of $\widetilde{G}$. Suppose that $\gcd(n,p)=1$. Then $\rho$ is not a dualizing involution, i.e., 
\[\pi^\rho\not\simeq\pi^\vee.\]
\end{theorem}

\begin{proof}
It follows from Lemma~\ref{pi rho f is same as pi f check rho} that $\pi^\rho\simeq\pi^\vee$ if and only if $\tr(\pi(\chi_{C^{\rho}}))=\tr(\pi(\chi_{C}))$ for all compact open subsets $C$ of $\widetilde{G}.$  
To establish the theorem, it suffices to show the existence of a compact open subset $C$ of $\widetilde{G}$ such that
\[
\tr(\pi(\chi_{C^\rho}))
\neq
\tr(\pi(\chi_C)).
\]
Since $(\pi,V)$ is admissible, there exists $m\geq 1$, and a compact open subgroup $D:=K_{m}^{*}$ such that  
\[
V^{D}\neq0,
\quad \text{and}  \quad 
\dim_{\mathbb{C}}(V^{D})<\infty.
\]
From Lemma~\ref{betatildefixesKm}, it follows that $\rho(D)=D$. Consider the projection $\pi(e_{D}): V\rightarrow V^{D}$ given by
\[
\pi(e_D)=\frac{1}{\mu(D)}\pi(\chi_D).
\]
Since $D$ is compact and open, we have $0<\mu(D)<\infty$ and 
\[
\tr(\pi(e_D))
=
\dim_{\mathbb{C}}(V^D)\neq 0.
\] 
Since $n\geq 3$, by Lemma~\ref{restriction-mu-n-of-central-character-injective}, there exists $\eta\in\mu_{n}$ such that $\omega_\pi(\eta^2)\neq1$, where $\omega_{\pi}$ is the central character of $\pi$. Let $C:=\eta D$. Since $\kappa(1_{K_{0}})=1$ and $\eta \neq 1$, it follows that $\eta\notin D$. Moreover,
\[
C^\rho=\rho(C)
=
\rho(\eta D)
=
\rho(D)\rho(\eta)
=
D\eta^{-1}
=
\eta^{-1}D
=
\eta^{-2}C.
\]
It is easy to see that $\pi(\chi_{C})=\omega_{\pi}(\eta)\pi(\chi_{D})$ and $\pi(\chi_{C^\rho})=\omega_\pi(\eta^{-2})\pi(\chi_C).
$
Indeed, 
\[ \pi(\chi_{C})= \int_{\widetilde{G}} \chi_{C}(g)\pi(g)dg = \int_{\eta D} \pi(g)dg = \int_D \pi(\eta g)dg = \omega_\pi(\eta)\pi(\chi_D),
\]
and 
\[
\pi(\chi_{C^{\rho}})= \int_{\widetilde{G}} \chi_{C^\rho}(g)\pi(g)dg = \int_{\eta^{-2}C} \pi(g)dg = \int_{C} \pi(\eta^{-2}g)dg = \omega_\pi(\eta^{-2})\pi(\chi_{C}).
\]
Since $\tr(\pi(\chi_{C}))\neq 0$ and $\omega_\pi(\eta^2)\neq 1$, 
we obtain
\[
\tr(\pi(\chi_{C^\rho}))
\neq
\tr(\pi(\chi_C)). \qedhere
\]
\end{proof}

\begin{corollary}\label{lifts of tau is not dualizing}
 Let $\pi$ be any irreducible admissible genuine representation of $\widetilde{G}$. Suppose that $\gcd(n,p)=1.$ For any $\alpha \in F^{\times}$, $\sigma_\alpha$ is not a dualizing involution, i.e.,
\[\pi^{\sigma_{\alpha}}\not\simeq \pi^{\vee}. \]
\end{corollary}

\begin{proof}  
Proceeding in a similar way as in  Theorem~\ref{Non-Dualizing lift}, choose $m\geq 1$ and a compact open subgroup $D:=K_{m}^{*}$ such that  \[ V^{D}\neq 0, \quad \text{and}  \quad \dim_{\mathbb{C}}(V^{D})<\infty. \] Since $\gcd(p,n)=1,$ by Hensel's lemma, we have $1+ \mathfrak{p} \subset (F^{\times})^{n}.$ Hence, $\det(x) \in (F^{\times})^{n}$ for every $x \in K_m.$ It follows that $$\sigma_{\alpha}(D) = (w_0, 1) D (w_0, 1)^{-1}.$$ Fix $\eta\in\mu_{n}$ such that $\omega_\pi(\eta^2)\neq1$ and let $C:=\eta D$. Then, we have $\sigma_\alpha(C)=(w_{0},1) \eta^{-2}C (w_{0},1)^{-1}$. Consequently, $\pi(\chi_{C^{\sigma_\alpha}}) = \pi(w_0, 1) \pi(\chi_{\eta^{-2}C}) \pi(w_0, 1)^{-1}.$ This implies that $$\tr(\pi(\chi_{C^{\sigma_\alpha}}))=\tr(\pi(\chi_{\eta^{-2}C})) =\omega_\pi(\eta^{-2}) \tr(\pi(\chi_C)).$$ 
It follows that $\tr(\pi(\chi_{C^{\sigma_\alpha}}))\neq \tr(\pi(\chi_C)).$ Hence the result.
\end{proof}

\subsubsection{Comparison of central characters (proof~2):} Assume that $r\geq 2$ and $n\geq 3$. Let $\rho$ be any lift of $\theta$ to $\widetilde{G}$. For an irreducible admissible  representation $\pi$ of $\widetilde{G},$ let $\omega_{\pi^{\rho}}$ and $\omega_{\pi^{\vee}}$ be the central characters of $\pi^{\rho}$ and $\pi^{\vee}$ respectively. 
\begin{theorem}[Non-dualizing lifts]\label{thm:non-dualizing-lifts}
Let $(\pi,V)$ be an irreducible admissible genuine representation of $\widetilde{G}$. Then $\rho$ is not a dualizing involution, i.e., 
\[\pi^{\rho}\not\simeq \pi^{\vee}.
\]
\end{theorem}
\begin{proof} Suppose that $\pi^{\rho}\simeq \pi^{\vee}$. For $\eta\in \mu_{n}$, and $v\in V$, we have 
\begin{equation*}
    \omega_{\pi^{\rho}}(\eta)v= \pi^\rho(\eta)v
    =\pi\left(\rho\left(\eta^{-1}\right)\right)v
    =\pi(\eta)v
    =\omega_{\pi}(\eta)v.
\end{equation*}
It follows that  
$$\omega_{\pi}(\eta)= \omega_{\pi^{\rho}}(\eta)=\omega_{\pi^{\vee}}(\eta)=\omega_{\pi}(\eta^{-1}).$$
Applying Lemma~\ref{restriction-mu-n-of-central-character-injective}, we get a contradiction. 
\end{proof}

\begin{corollary}\label{corollary:non-dualizing-lifts} Let $\pi$ be an irreducible admissible genuine representation of $\widetilde{G}$. For any $\alpha \in F^{\times}$, $\sigma_\alpha$ is not a dualizing involution, i.e., 
\[\pi^{\sigma_\alpha}\not\simeq \pi^{\vee}.
\]
\end{corollary}

\begin{proof} Suppose $\pi^{\sigma_\alpha}\simeq \pi^{\vee}$. For $\eta\in \mu_{n}$, it is easy to see that $\pi^{\sigma_{\alpha}}(\eta)=\pi^{\rho}(\eta)$, and hence we have 
$$\omega_{\pi^{\sigma_{\alpha}}}(\eta)=\omega_{\pi^{\rho}}(\eta).$$ The result follows. 
\end{proof}

\begin{remark} Let $r\geq 2$ and $n\geq 3$. Suppose that $\rho$ is any arbitrary involutory anti-automorphism (not necessarily a lift of $\theta$) on $\widetilde{G}$ satisfying $\rho(\eta)=\eta^{-1}$ for $\eta\in \mu_{n}$. Proceeding as in Theorem~\ref{thm:non-dualizing-lifts}, it can be shown by a simple comparison of central characters that such a $\rho$ can never be a duality involution. 
\end{remark}

\subsection{Dualizing lifts of $\tau$}

Let $r\geq 2$ be a positive integer and $G=\GL(r,F)$. Let $\widetilde{G}$ be the $2$-fold Kazhdan-Patterson cover of $G$ with respect to the $2$-cocycle $c\in \coh^{2}(G, \mu_{2})$. For $g\in G$, let $\theta(g)=~^{\top}g$, $\beta(g)=~^{\top}g^{-1}$ and $\tau(g)=w_{0} {}^{\top}g w_{0}$. Fix a lift $\rho$ of $\theta$ to $\widetilde{G}$ (for example, we can take $\rho$ as defined in Lemma~\ref{definition of the involution rho}). Let $G_{reg}$ be the subset of regular semisimple elements in $G$ and $\widetilde{G}_{reg}= pr^{-1}(G_{reg})$ be the subset of regular semisimple elements in $\widetilde{G}$. In this section, we prove the main results about dualizing involutions on $\widetilde{G}$. Before we proceed further, we record a few lemmas.

\begin{lemma} Let $g\in G_{reg}$. Then $F[g]$ is self-conjugate under the transpose operator, i.e., there exists $x\in G$ such that 
\[^{\top}F[g]=xF[g]x^{-1}. \]
\end{lemma}

\begin{proof} By \cite[Lemma~2]{Tupan2006}, there exists $x\in G$ such that $^{\top}x=x$ and $^{\top}g=xgx^{-1}$. For $t=f(g)\in F[g]$, we have 
\[^{\top}t= ~^{\top}f(g)=f(^{\top}g)=f(xgx^{-1})=xf(g)x^{-1}=xtx^{-1}.\] 
Hence the result. \qedhere

\end{proof}

\begin{remark}\label{xz is a symmetric conjugator of g} Let $z\in F[g]^\times$. Since $F[g]$ is self-conjugate under transpose, we get
\begin{align*}
^{\top}(xz)
&= ~^{\top}z^{\top}x=xzx^{-1}x=xz.
\end{align*}
Using $F[g]^{\times}=Z_{G}(g)$, we have $zg=gz$, and $(xz)g(xz)^{-1}= ~^{\top}g$. \\
\end{remark}

For $u,v\in F[g]^\times$, choose arbitrary lifts $\tilde{u},\tilde{v}\in\widetilde{G}$ and define
\begin{equation}\label{definition of B}
B(u,v):=[\tilde{u},\tilde{v}]
=\tilde{u}\tilde{v}\tilde{u}^{-1}\tilde{v}^{-1}.
\end{equation}
Since $F[g]^{\times}$ is abelian, the projection  $pr([\tilde{u},\tilde{v}])$ to $G$ is trivial. It follows that 
$B(u,v)\in\mu_2$. 
Clearly $B$ is independent of the chosen lifts since the kernel is
central, i.e., $\ker(pr) = \mu_2$. A simple computation yields $B(u,v)=c(u,v) c(v,u).$ Therefore, $$B:F[g]^\times \times F[g]^\times \longrightarrow \mu_2; \quad (u, v) \longmapsto B(u, v)$$ is a well defined map.

\begin{lemma}\label{properties of the map B} The map $B$ satisfies the following:
\begin{enumerate}
\item[\emph{1)}] $B(u, u) = 1$. 
\smallskip
\item[\emph{2)}] $B$ is symmetric.
\smallskip
\item[\emph{3)}] $B$ is bimultiplicative. 
\end{enumerate}
\end{lemma}

\begin{proof} Taking $u=v$, we have $B(u,u)=1$.  Clearly, $B$ is symmetric.  Taking $\tilde{u}_1  \tilde{u}_2$ as a lift of $u_1 u_2$, we have
\begin{align*}
    B(u_1u_2,v) &= \tilde{u}_1  \tilde{u}_2 \tilde v \tilde{u}_2^{-1} (\tilde{v}^{-1} \tilde{v})\tilde{u}_1^{-1} \tilde{v}^{-1} \\
    &= B(u_1,v) B(u_2,v).
\end{align*}
Since $B$ is symmetric, it follows that $B(v, u_1u_2)=B(v,u_1)B(v,u_2).$ 
\end{proof}

Let $g\in G_{reg}$. Let $x\in G$ such that $~^{\top}x=x$ and $xgx^{-1}= ~^{\top}g$. We call such an element $x$ a \textit{symmetric conjugator} of $g$. For $t\in F[g]^{\times}$, choose 
 $\tilde t\in pr^{-1}(t)$.  Since $\rho$ is a lift of $\theta$, we have $pr(\rho(\tilde{t}))= ~^{\top}t.$
On the other hand, for $\tilde{x}\in pr^{-1}(x)$, we get $$pr(\tilde{x}\tilde{t}\tilde{x}^{-1})
=xtx^{-1}= ~^{\top}t.$$
Thus $\rho(\tilde{t})$ and $\tilde{x}\tilde{t}\tilde{x}^{-1}$ differ by a unique element in $\mu_{2}$.

\begin{definition}\label{def:q-x} Let $g\in G_{reg}$ and $t\in F[g]^{\times}$. For $x$ a symmetric conjugator of $g$, we define
\begin{equation*}
q_x:F[g]^\times\longrightarrow\mu_2
\end{equation*}
as the unique element in $\mu_{2}$ satisfying \begin{equation*}
\rho(\tilde{t})
=q_x(t)\tilde{x}\tilde{t}\tilde{x}^{-1}.
\end{equation*}
\end{definition}

\begin{remark} Note that $q_{x}$ is independent of the chosen lifts. Indeed, replacing $\tilde{x}$ by $-\tilde{x}$ does not change the conjugation, while
replacing $\tilde{t}$ by $-\tilde{t}$ multiplies both sides by $-1$, 
since $\rho$ fixes $\mu_{2}$ pointwise.
\end{remark}

\begin{lemma}\label{lem:change-conjugator} For $z,t\in F[g]^{\times}$, we have
\[q_{xz}(t)=q_x(t)B(z,t).\]
\end{lemma}

\begin{proof} Consider the lift $\tilde{x}\tilde{z}$ of $xz$. Applying  \eqref{definition of B}, we have
\[
(\tilde{x}\tilde{z})\tilde{t}
(\tilde{x}\tilde{z})^{-1}
=\tilde{x}(\tilde{z}\tilde{t}\tilde{z}^{-1})\tilde{x}^{-1}
=B(z,t)\tilde{x}\tilde{t}\tilde{x}^{-1}.
\]
Comparing
\[
\rho(\tilde{t})
=q_x(t)\tilde{x}\tilde{t}\tilde{x}^{-1}
=q_{xz}(t)(\tilde{x}\tilde{z})
\tilde{t}(\tilde{x}\tilde{z})^{-1}
\]
gives
\begin{equation*}
q_x(t)=q_{xz}(t)B(z,t).
\end{equation*}
Since $B(z,t)^2=1$, the result follows.
\end{proof}

\begin{lemma}\label{lem:quadratic-q} For $t,u\in F[g]^{\times}$, we have
\begin{equation*}
q_x(tu)=q_x(t)q_x(u)B(t,u).
\end{equation*}
In particular, $q_{x}(t^{2})=1$.
\end{lemma}

\begin{proof}
Consider a lift $\tilde{t}\tilde{u}$ of $tu$. Using \eqref{definition of B} and $\rho$ is an anti-automorphism,  we have 
\begin{align*}
\rho(\tilde{t}\tilde{u})
&=\rho(\tilde{u})\rho(\tilde{t})\\
&=q_{x}(u)q_{x}(t)\tilde{x}
(\tilde{u}\tilde{t})\tilde{x}^{-1}\\
&= q_{x}(t)q_{x}(u)B(t,u)
\tilde{x}(\tilde{t}\tilde{u})\tilde{x}^{-1}\\
&= q_{x}(tu)\tilde{x}(\tilde{t}\tilde{u})\tilde{x}^{-1}.
\end{align*}
Thus, 
\[q_x(tu)=q_x(t)q_x(u)B(t,u).\]
Taking $u=t$ and using $B(t,t)=1$ gives $q_{x}(t^{2})=1$. 
\end{proof}

Let $u,v \in F[g]^{\times}$. In what follows, using Lemma~\ref{identifying F[g]^x as a product of Fi's}, we identify $u=(u_{1},\ldots,u_{s})$ and $v=(v_{1},\ldots,v_{s})$ for some $u_{i}, v_{i}\in F_{i}^{\times}$ for $i=1, 2, \dots, s$.

\begin{lemma}\label{description of B using norm} For $u, v \in F[g]^{\times}$, we have 
\[B(u,v) =\prod_{i=1}^{s}\bigl\langle u_iN(u),v_i\bigr\rangle_{F_i}.\]
\end{lemma}

\begin{proof} The embedding $\prod\limits_{i=1}^{s} F_i\hookrightarrow \M(r,F)$ is the multiplicity one
embedding: the $F_i$-block has $F$-dimension
$d_i=[F_i:F]$ and $d_1+\cdots+d_s=r.$  Specializing \cite[Proposition~4.2(5a), formula~(4.10)]{Zou2025}
to the standard Kazhdan--Patterson double cover, namely
$n=2$, $\boldsymbol c=0$, $\boldsymbol d=1$, gives
\begin{equation*}
B(u,v)
=\prod_{i=1}^s \langle u_i,v_i\rangle_{F_i}
\bigl\langle\det(u),\det(v)\bigr\rangle_{F}.
\end{equation*}
By Lemma~\ref{lem:det-equals-norm}, we get
\begin{equation*}
B(u,v)
=\prod_{i=1}^s \langle u_i,v_i\rangle_{F_i}
\bigl\langle N(u),N(v)\bigr\rangle_{F}.
\end{equation*}
Since $N(v)=\prod\limits_{i=1}^{s} N_{F_i/F}(v_i),$ using
bimultiplicativity of the Hilbert symbol and Lemma~\ref{relation between Hilbert symbol on K and F}, we have
\begin{align*}
\bigl\langle N(u),N(v)\bigr\rangle_{F}
&=\prod_{i=1}^{s} \bigl\langle N(u),N_{F_i/F}(v_i)\bigr\rangle_{F}=\prod_{i=1}^{s} \bigl\langle N(u),v_i\bigr\rangle_{F_i}.
\end{align*}
Thus
\begin{align*}
B(u,v)
&=\prod_{i=1}^s \langle u_i,v_i\rangle_{F_i}
\prod_{i=1}^s \bigl\langle N(u),v_i\bigr\rangle_{F_i}=\prod_{i=1}^s\bigl\langle u_iN(u),v_i\bigr\rangle_{F_i}. \qedhere
\end{align*}
\end{proof}
\vspace{0.5 cm}

\noindent Consider the radical of the bimultiplicative map $B$ given by
\begin{equation}\label{Rad(B)}
\Rad(B):=\left\{u\in F[g]^\times:B(u,v)=1, 
\forall v\in F[g]^\times\right\}.
\end{equation}
Since $B(v,u)=B(u,v)$, the left and right radicals
coincide.

\begin{lemma}\label{lem:image-B} For $a\in F[g]^{\times}$, let $\chi_{a}:F[g]^{\times}\longrightarrow \mu_{2}$ be the map
\[\chi_{a}(z)=B(z,a).\]
Then $\chi_a$ is a character and 
$$\operatorname{Im}(\chi_{a})=\begin{cases}
\{1\}, & a\in\Rad(B)\\
\mu_{2}, & a\notin\Rad(B).
\end{cases}$$
\end{lemma}

\begin{proof}
Since $B$ is bimultiplicative, it is clear that $\chi_{a}$ is a character.  If
$a\in\Rad(B)$, using $B(a,z)=1$ for every $z \in F[g]^{\times}$ and $B$ is symmetric, we have $\operatorname{Im}(\chi_a)=\{1\}.$ Suppose that $a\notin\Rad(B)$. There exists $z\in F[g]^{\times}$ such that $B(a,z)\neq 1$. Since $B$ takes values in $\mu_{2}$, we have $B(z,a)=B(a,z)=-1$. It follows that $\chi_{a}$ is surjective. 
\end{proof}

Having an explicit description of the radical $\Rad(B)$ is important in showing that for $\tilde{g}\in \widetilde{G}_{reg}$, we have $\rho(\tilde{g})$ is a conjugate of $\tilde{g}$. The fact that $\rho(\tilde{g})$  and $\tilde{g}$ are in the same conjugacy class in $\widetilde{G}$ is crucial in establishing the dualizing property of $\rho$.

\subsubsection{Even rank case:} \label{Even rank case}
Let $r=2t$ be an even integer and $G=\GL(r,F)$. Let $\tau$ be the standard involution on $G$. In this section, we show that any lift of $\tau$ to $\widetilde{G}$ is a dualizing involution. \\

\begin{lemma}\label{prop:even-radical} If $r=2t$ is even, then 
\[\Rad(B)=(F[g]^\times)^{2}.\]
\end{lemma}

\begin{proof} We show that $\Rad(B)\subseteq (F[g]^\times)^{2}$ and vice versa. Let $u=(u_1,\ldots,u_s)\in\Rad(B)$. Fix $1\leq j\leq s$ and let  $v_j\in F_j^\times$. Consider $v=(1,\ldots,1,v_j,1,\ldots,1) \in F[g]^\times$. Applying Lemma~\ref{description of B using norm}, we have 
\begin{equation*}
\bigl\langle u_jN(u),v_j\bigr\rangle_{F_j}=B(u,v)=1.
\end{equation*}
 Using non degeneracy of the Hilbert
symbol, it follows that 
\begin{equation}\label{ujN(u) is a square}
u_{j}N(u)\in (F_{j}^{\times})^{2}, \quad \text{for each} \quad 1\leq j \leq s.\end{equation}
Let $[a]$ denote the class of $a$ in the relevant square class group.  Using \eqref{ujN(u) is a square}, we get  
\begin{equation*}
[u_j]=[N(u)]
\quad \text{in }F_j^\times/({F_j}^\times)^{2}.
\end{equation*}
By multiplicative property of the norm, we have $N_{K/F}(b a^2)=N_{K/F}(b)N_{K/F}(a)^2$ and it follows that the norm induces a map on the square classes. Applying 
$N_{F_j/F}$ to the preceding equality yields
\begin{equation*}
[N_{F_j/F}(u_j)]
=[N_{F_j/F}(N(u))]=[N(u)]^{d_j}.
\end{equation*}
Multiplying over all $j$ gives
\begin{align*}
[N(u)]
&=\left[\prod_{j=1}^s N_{F_j/F}(u_j)\right]=[N(u)]^{d_1+\cdots+d_s}=[N(u)]^r.
\end{align*}
Since $r$ is even, $N(u)\in (F^{\times})^{2}$. Using \eqref{ujN(u) is a square}, it follows that $u_{j}\in (F_{j}^{\times})^{2}$ for each $j$. Thus \[u\in (F[g]^\times)^2.\] 
Conversely, suppose that $y=b^2$ for some $b\in F[g]^\times$.  By
bimultiplicativity of $B$,
\begin{equation*}
B(y,v)=B(b^2,v)=B(b,v)^2=1
\end{equation*}
for every $v\in F[g]^\times$.  Thus
$y\in\Rad(B)$. Hence the result. 
\end{proof}

\begin{theorem}\label{thm:even-conjugacy} For $\tilde{g}\in \widetilde{G}_{reg}$, there exists $\tilde{y}\in\widetilde{G}$ such that 
\[\rho(\tilde{g})=\tilde{y}\tilde{g}\tilde{y}^{-1}.
\]
\end{theorem}

\begin{proof} Write $\tilde{g}=(g,\epsilon).$ Let $x\in G$ such that $^{\top}x=x$ and $^{\top}g=xgx^{-1}$.
We consider the cases when $g\notin \Rad(B)$ and $g\in\Rad(B)$. In both the cases, we show that $\rho(\tilde{g})$ and $\tilde{g}$ are in the same conjugacy class in $\widetilde{G}$. Suppose $g\notin\Rad(B).$ From Lemma~\ref{lem:image-B}, it follows that the character
\begin{equation*}
z\longmapsto B(z,g)
\end{equation*}
of $F[g]^{\times}$ has image $\mu_2$. Since $q_{x}(g)\in\mu_2$, there exists
$z\in F[g]^\times$ such that $B(z,g)=q_x(g)$. By Lemma~\ref{lem:change-conjugator},
\[
q_{xz}(g)=q_x(g)B(z,g)=q_x(g)^2=1.
\]
Since $^{\top}(xz)=(xz)$, and $^{\top}g=(xz)g(xz)^{-1}$ (see Remark~\ref{xz is a symmetric conjugator of g}), we have 
\[
\rho(\tilde{g})
=(\widetilde{xz})\tilde{g}(\widetilde{xz})^{-1}.
\]
Suppose that $g\in\Rad(B).$ By Proposition~\ref{prop:even-radical}, there exists $b\in F[g]^\times$ such that $g=b^2$, and Lemma~\ref{lem:quadratic-q}, gives
\[
q_{x}(g)=q_{x}(b^{2})=1.
\]
Consequently,
\[
\rho(\tilde{g})
=\tilde{x}\tilde{g}\tilde{x}^{-1}.  \qedhere
\] 
\end{proof}

\begin{corollary}\label{conjugacy statement wrt sigma} For $\tilde{g}\in \widetilde{G}_{reg}$, there exists $\tilde{u}\in\widetilde{G}$ such that 
\[\sigma(\tilde{g})=\tilde{u}\tilde{g}\tilde{u}^{-1}.
\]
\end{corollary}

\begin{proof} By Theorem~\ref{thm:even-conjugacy}, there exists $\tilde{y}\in\widetilde{G}$ such that $\rho(\tilde{g})=\tilde{y}\tilde{g}\tilde{y}^{-1}.$ Taking $\tilde{u} = (w_{0}, 1) \tilde{y}$, we get 
$$\sigma(\tilde{g})= (w_{0}, 1) \rho(\tilde{g}) (w_{0}, 1)^{-1} = \tilde{u}\tilde{g}\tilde{u}^{-1}.$$
Hence the result.
\end{proof}

\begin{corollary}\label{conjugacy statement wrt sigma alpha} Let $\tilde{g}\in \widetilde{G}_{reg}$ and $\alpha\in F^{\times}$. Let $\sigma_{\alpha}$ be the lift of $\tau$ defined in Lemma~\ref{other lifts of tau}. There exists $\tilde{u}\in\widetilde{G}$ such that 
\[\sigma_\alpha(\tilde{g})=\tilde{u}\tilde{g}\tilde{u}^{-1}.
\]
 \end{corollary}

\begin{proof} 
If $g\in \Rad(B),$ using Lemma~\ref{prop:even-radical} we have $\sigma_{\alpha}(\tilde{g})=\sigma(\tilde{g})$. The result follows from Corollary~\ref{conjugacy statement wrt sigma}. We may assume that $g\notin \Rad(B)$  and $\langle \alpha, \triangle(\tilde{g}) \rangle_{F} =-1$. Using Lemma~\ref{lem:image-B}, we can choose $z \in F[g]^{\times}$ such that $B(z, g) = \langle \alpha, \triangle(\tilde{g}) \rangle_{F}.$ Write $y =\langle \alpha, \triangle(\tilde{g}) \rangle_{F}\tilde{g}.$ By Corollary~\ref{conjugacy statement wrt sigma}, there exists $\tilde{v} \in \widetilde{G}$ such that $\sigma(y) = \tilde{v}y \tilde{v}^{-1}$. Thus $$\sigma_\alpha(\tilde{g})= \sigma(y)= \tilde{v}y \tilde{v}^{-1}.$$
Since $zg = gz,$ a simple computation gives $$B(z, g) = c(z, g) c(zg,z^{-1}) c(z, z^{-1}).$$ Hence $y = \tilde{z} \tilde{g} \tilde{z}^{-1},$ and the result follows.
\end{proof}

\begin{theorem}[Dualizing lifts]\label{main theorem 2} Let $\pi$ be an irreducible admissible genuine representation of $\widetilde{G}$ and $\rho$ be a lift of $\theta$ to $\widetilde{G}$. Then 
\[\pi^{\rho}\simeq \pi^{\vee}.\]
\end{theorem}

\begin{proof} To establish that $\rho$ is a duality involution, it is enough to show that 
\[\Theta_{\pi^{\vee}}(f)=\tr(\pi^\vee(f))=\tr(\pi^{{\rho}}(f))=\Theta_{\pi^{\rho}}(f)\] 
 for all $f \in  C_c^\infty(\widetilde G)$.
Using Lemma~\ref{pi rho f is same as pi f check rho}, we have
\begin{equation}\label{trace of pi rho f}
    \tr(\pi^{{\rho}}(f))=\tr(\pi((f^\vee)^{\rho}))=\int_{\widetilde{G}_{reg}} f(\rho(g^{-1}))\Theta_{\pi}(g)dg.
\end{equation}
Applying the change of variable $g\mapsto\rho(g^{-1})$ in \eqref{trace of pi rho f} and using $\rho(g^{-1})$ is conjugate to $g^{-1}$ (see Theorem~\ref{thm:even-conjugacy}), we have 
\[\tr(\pi^{\rho}(f))= 
\int_{\widetilde{G}_{reg}} f(g)\Theta_{\pi}(g^{-1})dg.\]
By part 3) of Lemma~\ref{pi rho f is same as pi f check rho},
it follows that
\[\tr(\pi^{\rho}(f)) 
=\int_{\widetilde{G}_{reg}} f(g)\Theta_{\pi^\vee}(g)dg = \tr({\pi^\vee}(f)).\qedhere\] 
\end{proof}

\begin{corollary} Let $\alpha \in F^{\times}$ and $\sigma_{\alpha}$ be any lift of $\tau$. Let $\pi$ be any irreducible admissible genuine representation of $\widetilde{G}$. Then 
\[\pi^{\sigma_{\alpha}}\simeq \pi^{\vee}. \]
\end{corollary}

\begin{proof} Using Corollary~\ref{conjugacy statement wrt sigma alpha}, and repeating the arguments in Theorem~\ref{main theorem 2} with $\rho$ replaced by $\sigma_{\alpha}$, it follows that $\sigma_{\alpha}$ is dualizing. 
\end{proof}

\subsubsection{Odd rank case:} \label{Odd rank case}

Let $r=2t+1$ be an odd integer. Let $G=\GL(r,F)$ and $\widetilde{G}$ be the Kazhdan-Patterson cover of $G$ with respect to the $2$-cocycle $c$. Let $\theta(g)=~^{\top}g$ and $\rho$ be any lift of $\theta$ to $\widetilde{G}$. For $\lambda\in F^{\times}$, we let $z_{\lambda}=\lambda1_{G}$ and write $\tilde{z}_{\lambda}=(z_{\lambda}, \epsilon) \in \widetilde{G}$ for any lift of $z_{\lambda}$. In the odd rank case, we show that there is a unique dualizing lift of the standard involution.

\begin{lemma}\label{definition of chirho} Let $\chi_{\rho}: F^{\times} \rightarrow \mu_{2}$ be the map defined by:
\[\chi_{\rho}(\lambda)= \rho(\tilde{z}_{\lambda})\tilde{z}_{\lambda}^{-1}. \]
Then $\chi_{\rho}$ is a well defined character of $F^{\times}$. 
\end{lemma}

\begin{proof}  We have 
\begin{align*}
\rho(\tilde{z}_{\lambda})\tilde{z}_{\lambda}^{-1} &= \epsilon^{-1}\rho(z_{\lambda}, 1)(z_{\lambda}^{-1}, \epsilon^{-1}c(z_{\lambda}, z_{\lambda}^{-1})^{-1})\\
&= \epsilon^{-2}\rho(z_{\lambda}, 1)(z_{\lambda}^{-1}, c(z_{\lambda}, z_{\lambda}^{-1}))\\
&= \rho(z_{\lambda, 1})(z_{\lambda}^{-1}, c(z_{\lambda}, z_{\lambda}^{-1}))\\
&= (z_{\lambda}, *)(z_{\lambda}^{-1}, c(z_{\lambda}, z_{\lambda}^{-1}))\\
&=(1_{G}, \bullet).
\end{align*}
Thus, $\chi_{\rho}$ is a well defined map and it is independent of the chosen lifts. For $\lambda_{1}, \lambda_{2} \in F^{\times},$ we have 
\begin{equation}\label{to conclude chi rho is a character}
\rho(\tilde{z}_{{\lambda}_{1}} \tilde{z}_{{\lambda}_{2}}) = \chi_{\rho}(\lambda_{1}) \chi_{\rho}(\lambda_{2}) \tilde{z}_{{\lambda}_{2}} \tilde{z}_{{\lambda}_{1}}\quad \text{and} \quad \rho(\tilde{z}_{{\lambda}_{1} {\lambda}_{2}}) = \chi_{\rho}(\lambda_{1} \lambda_{2}) \tilde{z}_{{\lambda}_{1} {\lambda}_{2}}.
\end{equation}
Since $ \tilde{z}_{{\lambda}_{1}} \tilde{z}_{{\lambda}_{2}} = \tilde{z}_{{\lambda}_{2}} \tilde{z}_{{\lambda}_{1}},$ using $\tilde{z}_{{\lambda}_{1}} \tilde{z}_{{\lambda}_{2}}$ as a lift of $z_{{\lambda}_{1} {\lambda}_{2}}$ in \eqref{to conclude chi rho is a character}, it follows that 
$$\chi_{\rho}(\lambda_{1} \lambda_{2}) = \chi_{\rho}(\lambda_{1})\chi_{\rho}(\lambda_{2}).$$ Hence the result.
\end{proof}

\begin{lemma}\label{definition of rhohash} Let $\rho^{\#}: \widetilde{G} \to \widetilde{G}$ be the map given by:
\begin{equation}\label{definition of rho sharp}
\rho^{\#}(g, \epsilon) = \chi_{\rho}(\det(g))\rho(g, \epsilon).
\end{equation}
Then 
\begin{enumerate}[label=\emph{\arabic*)}]
\item $\rho^{\#}$ is an anti-automorphism on $\widetilde{G}$ and $(\rho^{\#})^2 =1.$
\smallskip

\item $\rho^{\#}$ is a lift of $\theta$.
\end{enumerate}
\end{lemma}
The proof follows from Proposition~\ref{description of the set of liftings}.

\begin{lemma} \label{rho hash fixes zlambda} For $\lambda\in F^{\times}$, we have 
$$\rho^{\#}(\tilde{z}_{\lambda}) = \tilde{z}_{\lambda}.$$
\end{lemma}

\begin{proof} Using $r=2t+1$ is odd, and $\chi_{\rho}(\lambda)\in \mu_{2}$, we get $\rho^{\#}(\tilde{z}_{\lambda})= \tilde{z}_{\lambda}$. Indeed,
\begin{align*} \rho^{\#}(\tilde{z}_{\lambda}) &=\chi_\rho(\det(z_{\lambda}))\rho(\tilde{z}_{\lambda}) \\
&=\chi_\rho(\lambda^r) \chi_\rho(\lambda)\tilde{z}_{\lambda}\\
&=\chi_\rho(\lambda)^{r+1} \tilde{z}_{\lambda}\\ 
&=\tilde{z}_{\lambda}.\qedhere
\end{align*}
\end{proof}

\begin{lemma}\label{prop:odd-radical}
For $g \in G_{reg},$ we have
\[\Rad(B)= F^{\times}(F[g]^{\times})^{2}.\]
\end{lemma}

\begin{proof} 
Let $u=(u_1,\ldots,u_s)\in\operatorname{Rad}(B)$. Proceeding in a similar way as in Lemma~\ref{prop:even-radical}, we get
\begin{equation*}
u_iN(u)\in (F_i^{\times})^{2}, \quad \text{for every} \quad 1\leq i \leq s.
\end{equation*}
For each $i$, choose $b_i\in F_i^\times$ such that $u_i=N(u)b_i^2$. Since $b:=(b_1,\ldots,b_s)\in F[g]^{\times}$, it follows that $u\in F^{\times}(F[g]^{\times})^{2}$.
Conversely, suppose that $y=\lambda a^2$ for $\lambda\in F^{\times}$ and $a\in F[g]^{\times}$. Since $\dim_F F[g]=r$, we have $N(\lambda)=\lambda^r.$ Using the multiplicativity of the norm, we have
\begin{align*}
y N(y)
&=\lambda a^{2}\lambda^{r}N(a)^{2}
= \bigl(\lambda^{r+1} N(a)^{2} \bigr) a^{2}.
\end{align*}
Since $r+1$ is even, it is easy to see that $y N(y)$ is a square in $F[g]^{\times}$. Writing $y=(y_1,\ldots,y_s),$ it follows that $y_{i}N(y)$ is a square in $F_{i}^{\times}$ for each $i$.  Using Lemma~\ref{description of B using norm}, we have  $$B(y,v)=1,$$
for every $v\in F[g]^\times$.  Thus $y\in\operatorname{Rad}(B)$.
\end{proof}

\begin{theorem}\label{thm:odd-conjugacy} For $\tilde{g}\in \widetilde{G}_{reg}$,  there exists $\tilde{y}\in\widetilde{G}$ such that 
\[\rho^{\#}(\tilde{g})=\tilde{y}\tilde{g}\tilde{y}^{-1}.
\]
\end{theorem}

\begin{proof} We write  $\tilde{g}=(g,\epsilon)$. Let $x\in G$ such that $^{\top}x=x$ and $^{\top}g=xgx^{-1}$. We consider the cases when $g\notin \Rad(B)$ and $g\in\Rad(B)$. In both the cases, we show that $\rho(\tilde{g})$ and $\tilde{g}$ are in the same conjugacy class in $\widetilde{G}$. Suppose that $g\notin\Rad(B)$. 
By Lemma~\ref{lem:image-B}, it follows that the character
\begin{equation*}
z\longmapsto B(z,g)
\end{equation*}
of $F[g]^{\times}$ has image $\mu_{2}$. We choose $z\in F[g]^{\times}$ such that $B(z,g)q_{x}(g)=\chi_{\rho}(\det(g))$. Using Lemma~\ref{lem:change-conjugator}, we have
\[
q_{xz}(g)=q_{x}(g)B(z,g)=\chi_{\rho}(\det(g)).
\]
Since $^{\top}(xz)=(xz)$, and $^{\top}g=(xz)g(xz)^{-1}$, it is easy to see that 
\[
\rho^{\#}(\tilde{g})
=(\widetilde{xz})\tilde{g}(\widetilde{xz})^{-1}.
\]
Indeed, 
\begin{align*}
\rho^{\#}(\tilde{g}) & = \chi_{\rho}(\det(g))\rho(\tilde{g}) \\
&= \chi_{\rho}(\det(g)) q_{xz}(g) (\widetilde{xz})\tilde{g}(\widetilde{xz})^{-1} \\
&= \chi_{\rho}(\det(g))^{2} (\widetilde{xz})\tilde{g}(\widetilde{xz})^{-1} \\
&= (\widetilde{xz})\tilde{g}(\widetilde{xz})^{-1}.
\end{align*}
Suppose that $g\in \Rad(B).$ By Proposition~\ref{prop:odd-radical}, we can write $g=z_{\lambda}b^{2},$ for some $\lambda\in F^{\times}$ and $b \in F[g]^{\times}.$ Using Lemma~\ref{lem:quadratic-q}, we have
\begin{align*}
\rho^{\#}(\tilde{g})
&=\chi_{\rho}(\det(z_{\lambda} b^{2}))q_{x}(g)\tilde{x}\tilde{g}\tilde{x}^{-1}\\
&=\chi_{\rho}(\lambda)q_{x}(z_{\lambda})q_{x}(b^{2})B(z_{\lambda},b^{2})\tilde{x}\tilde{g}\tilde{x}^{-1}\\
&=\chi_{\rho}(\lambda)q_{x}(z_{\lambda})\tilde{x}\tilde{g}\tilde{x}^{-1}\\
&=\rho(\tilde{z}_{\lambda})\tilde{z}_{\lambda}^{-1}q_{x}(z_{\lambda})\tilde{x}\tilde{g}\tilde{x}^{-1}\\
&= (q_{x}(z_{\lambda})\tilde{x}\tilde{z}_{\lambda}\tilde{x}^{-1}) \tilde{z}_{\lambda}^{-1} (q_{x}(z_{\lambda})\tilde{x}\tilde{g}\tilde{x}^{-1}) \\
& = q_{x}(z_{\lambda})^{2}(\tilde{x}\tilde{z}_{\lambda}\tilde{x}^{-1}\tilde{z}_{\lambda}^{-1}) \tilde{x}\tilde{g}\tilde{x}^{-1}.
\end{align*}
Since $\tilde{z}_{\lambda}\in Z(\widetilde{G})$, we get $\rho^{\#}(\tilde{g})=\tilde{x}\tilde{g}\tilde{x}^{-1}.$ Hence the result.  
\end{proof}

\begin{theorem}[Dualizing lifts]\label{main theorem odd case} Let $\pi$ be an irreducible admissible genuine representation of $\widetilde{G}$ and $\rho^{\#}$ be the lift of $\theta$ to $\widetilde{G}$. Then 
\[\pi^{\rho^{\#}}\simeq \pi^{\vee}.\]
\end{theorem}

\begin{proof} Using Theorem~\ref{thm:odd-conjugacy}, and repeating the arguments in Theorem~\ref{main theorem 2} with $\rho$ replaced by $\rho^{\#}$, it follows that $\rho^{\#}$ is a dualizing involution. 
\end{proof}

\begin{lemma}\label{chirho is trivial iffn} $\rho$ is a dualizing involution if and only if $\chi_\rho$ is trivial.
\end{lemma}

\begin{proof} If $\chi_{\rho}$ is trivial, then $\rho(\tilde{g})=\rho^{\#}(\tilde{g})$ for all $\tilde{g}\in \widetilde{G}$ and the result follows from Theorem~\ref{main theorem odd case}. For the converse, we show that if $\chi_{\rho}$ is non-trivial, then $\rho$ is not dualizing. Suppose $\chi_{\rho}\neq 1$ and $\rho$ is dualizing. There exists $\lambda \in F^{\times}$ such that $\chi_{\rho}(\lambda)=-1$. Thus, we have $\rho(\tilde{z}_{\lambda})= -\tilde{z}_{\lambda}.$ Consequently, for an irreducible admissible genuine representation $\pi$ of $\widetilde{G},$
$$\pi^{\rho}(\tilde{z}_{\lambda}) = \pi(\rho(\tilde{z}_{\lambda}^{-1}))= \pi(-\tilde{z}_{\lambda}^{-1}) = -\pi(\tilde{z}_{\lambda}^{-1}).$$ 
Since $\tilde{z}_{\lambda}\in Z(\widetilde{G})$, the preceding equality yields $\omega_{\pi^{\rho}}(\tilde{z}_{\lambda})= -\omega_{\pi}(\tilde{z}_{\lambda}^{-1})$. On the other hand, as $\pi^{\rho} \simeq \pi^{\vee}$, we get $\omega_{\pi^{\rho}}(\tilde{z}_{\lambda})= \omega_{\pi}(\tilde{z}_{\lambda}^{-1}).$ Hence, we have $-\omega_\pi(\tilde{z}_{\lambda}) = \omega_\pi(\tilde{z}_{\lambda}).$ This is a contradiction. 
\end{proof}

\begin{theorem}[Uniqueness]\label{Uniqueness of dualizing involution} Let $\pi$ be an irreducible admissible genuine representation of $\widetilde{G}$. Then there exists a unique lift $\rho$ of $\theta$ to $\widetilde{G}$ such that
\[\pi^{\rho}\simeq \pi^{\vee}.\]
\end{theorem}

\begin{proof} The existence of a lift $\rho$ of $\theta$ to $\widetilde{G}$ which is dualizing follows from Theorem~\ref{main theorem odd case}. Suppose that $\rho$ and $\rho'$ are dualizing involutions on $\widetilde{G}$. Applying Lemma~$\ref{chirho is trivial iffn}$, we have 
\begin{equation}\label{rho and rho' fix the center}
\rho(\tilde{z}_{\lambda})=\tilde{z}_{\lambda}=\rho'(\tilde{z}_{\lambda}), \text{ for all } \lambda \in F^{\times}.
\end{equation}
Since any character $\phi:G \to \mu_{2}$ is given by $x \mapsto \psi(\det(x) ),$ where  $\psi: F^{\times} \to \mu_{2}$ is a character, using Proposition~\ref{description of the set of liftings}, we have $$\rho'(\tilde{g})=\psi(\det(g))\rho(\tilde{g}), \text{ for all } \tilde{g} = (g, \epsilon) \in \widetilde{G}.$$
In particular, taking $\tilde{g} = \tilde{z}_{\lambda}$ in the preceding equality yields
\begin{equation}\label{equation to conclude the uniquness of the involution}
\rho'(\tilde{z}_{\lambda})=\psi(\lambda)\rho(\tilde{z}_{\lambda}).
\end{equation}
By \eqref{rho and rho' fix the center} and \eqref{equation to conclude the uniquness of the involution}, we have $\psi(\lambda)=1$ for all $\lambda\in F^{\times}.$ The result follows. 
\end{proof}

\begin{corollary}\label{Uniqueness of lift sigma} Let $\pi$ be an irreducible admissible genuine representation of $\widetilde{G}$.  Then there exists a unique lift $\sigma$ of $\tau$ to $\widetilde{G}$ such that
\[\pi^{\sigma}\simeq \pi^{\vee}.\]
\end{corollary}

\begin{proof}
By Theorem~\ref{Uniqueness of dualizing involution}, there exists a unique lift $\rho$ of $\theta$ to $\widetilde{G}$ which is dualizing. Define $\sigma: \widetilde{G} \to \widetilde{G}$ as follows: 
\[\sigma(g)=(w_0,1)\rho(g)(w_0,1)^{-1}.\]
Using Lemma~\ref{sigma-rho dualizing property}, we have that $\sigma$ is dualizing. Suppose there exists $\alpha \in F^\times$ such that $\sigma_{\alpha}\neq \sigma$ and ${\sigma_\alpha}$ is dualizing lift of $\tau$. Consider the map $\rho_{\alpha}: \widetilde{G} \to \widetilde{G}$ given by $$\rho_\alpha(g)=(w_0,1)^{-1}\sigma_{\alpha}(g)(w_0,1).$$ Clearly, $\rho_\alpha$ is a lift of $\theta$ to $\widetilde{G}$, and is a dualizing involution. However $\rho_\alpha \neq \rho$, contradicting the uniqueness of $\rho.$ Hence the result. \qedhere
\end{proof}

\section{Dualizing involutions on twisted Kazhdan-Patterson covers of $G$} \label{Twisted covers}

Throughout, we let $F$ be a non-Archimedean local field of characteristic zero and  $p$ be the residual degree of $F.$ Let $r, n \geq 2$ be  positive integers and $G=\GL(r,F)$. For $m\in \mathbb{Z}/n\mathbb{Z},$ let $c^{(m)}\in \mathrm{H}^{2}(G, \mu_{n})$ be the $2$-cocycle given in Section~\ref{Kazhdan-Patterson covers of GL(r)} and $\widetilde{G}^{(m)}$ be the central extension of $G$ with respect to the $2$-cocycle $c^{(m)}$. 
Let $\theta$ be the transpose anti-automorphism on $G$. Let $\tau$ be the standard involution on $G.$
In this section, we state and prove results about dualizing lifts of $\tau$ to $\widetilde{G}^{(m)}$. \\

Fix $m\in \mathbb{Z}/n\mathbb{Z}$ and denote $\Lambda:=c^{(m)}$. Let $\beta$ be the transpose inverse automorphism on $G$. Using \eqref{sigma beta and sigma cohologous main equation}, it follows that the cocycle $\Lambda^\beta$ is cohomologous to $\Lambda$ via the same map $\delta$. Thus, Lemma~\ref{to show beta tilde is an involution} holds. Moreover, the maps $\tilde{\beta}:\widetilde{G}^{(m)}\rightarrow \widetilde{G}^{(m)}$ and $\rho: \widetilde{G}^{(m)} \to \widetilde{G}^{(m)}$ defined by $$\tilde{\beta} (x, \epsilon) = (\beta(x), \epsilon \delta(x)) \quad \text{and} \quad \rho(x, \epsilon) = \tilde{\beta}((x, \epsilon)^{-1})$$ satisfy analogous properties as in Lemma~\ref{lift of the automorphism beta} and Lemma~\ref{definition of the involution rho}.
Let $\sigma: \widetilde{G}^{(m)} \to \widetilde{G}^{(m)}$ be the map given by 
\[\sigma(g)=(w_{0},1) \rho(g) (w_{0},1)^{-1}.\]
For $\alpha \in F^{\times},$ define $\sigma_\alpha: \widetilde{G}^{(m)} \to \widetilde{G}^{(m)}$ as follows: 
$$\sigma_{\alpha}(g)= \langle \alpha, \triangle(g)\rangle_{F} \,\sigma(g).$$ 
Using similar arguments as in  Lemma~\ref{definition of the involution sigma} and Lemma~\ref{other lifts of tau}, we can show that $\sigma$ and $\sigma_{\alpha}$ are lifts of $\tau$ to $\widetilde{G}^{(m)}.$ Moreover, by Proposition~\ref{description of the set of liftings}, it follows that for $\alpha \in F^{\times},$  $\sigma_\alpha$ are all possible lifts of $\tau$ to $\widetilde{G}^{(m)}.$ \\

To show that a lift $\rho$ of $\theta$ to $\widetilde{G}^{(m)}$ is either  non-dualizing or dualizing is equivalent to showing that a lift $\sigma_\alpha$ of $\tau$ to $\widetilde{G}^{(m)}$ is either non-dualizing or dualizing (see Lemma~\ref{sigma-rho dualizing property}). We now present the main results on the dualizing property of the lifts $\sigma_\alpha$.

\begin{theorem}[Non-Dualizing lifts] Suppose that $r\geq 2$ and $n\geq3$. Let $G=\GL(r,F)$ and let $\widetilde{G}^{(m)}$ be the $n$-fold Kazhdan-Patterson cover of $G$ with respect to $\Lambda.$ Let $(\pi, V)$ be an irreducible admissible genuine representation of $\widetilde{G}^{(m)}$.  For any $\alpha \in F^{\times}$, $\sigma_\alpha$ is not a dualizing involution, i.e.,
\[\pi^{\sigma_{\alpha}}\not\simeq \pi^{\vee}. \] 
\end{theorem}

Proceeding as in Corollary~\ref{corollary:non-dualizing-lifts}, we can establish the result. \\

In what follows, we fix $n=2$ and let  $\widetilde{G}^{(m)}$ denote the $2$-fold cover of $G$ with respect to $\Lambda \in \coh^{2}(G, \mu_{2})$. Let $g \in G_{reg}.$ For $u,v\in F[g]^\times$, choose arbitrary lifts $\tilde{u},\tilde{v}\in\widetilde{G}^{(m)}$ and define $B: F[g]^\times \times F[g]^\times \longrightarrow \mu_2$ as follows: 
$$B(u,v):=[\tilde{u},\tilde{v}]
=\tilde{u}\tilde{v}\tilde{u}^{-1}\tilde{v}^{-1}.$$

Then Lemma~\ref{properties of the map B}, Definition~\ref{def:q-x}, Lemma~\ref{lem:change-conjugator} and Lemma~\ref{lem:quadratic-q} continue to hold in the setup of $\widetilde{G}^{(m)}$.

\begin{lemma}\label{Rad for twisted covers}
For $g \in G_{reg},$ we have
$$\Rad(B) = \begin{cases}
     (F[g]^\times)^2, \quad \mbox{if $r$ is even} \\ 
     F^\times (F[g]^\times)^2, \quad \mbox{if $r$ is odd} .
 \end{cases}$$

\end{lemma}

\begin{proof}
Note that Lemma~\ref{description of B using norm} will hold in the context of $\widetilde{G}^{(m)}$ due to the following reasons:
   \begin{enumerate}[label=\roman*), itemsep=0.2cm]
        \item We will apply \cite[Proposition~4.2(5a), formula~(4.10)]{Zou2025} with $\boldsymbol c=m$, $\boldsymbol d=1$ and $n=2$. 

        \item For $u, v \in F[g]^\times,$ the factor $\langle \det_F(u),\det_F(v) \rangle_F^{2m+1}$ appearing in \cite[Proposition~4.2(5a), formula~(4.10)]{Zou2025} is equal to $\langle \det_F(u),\det_F(v) \rangle_F$ as $n=2.$
  \end{enumerate}
Thus, the result follows repeating the arguments used in  Lemma~\ref{prop:even-radical} and  Lemma~\ref{prop:odd-radical}.
\end{proof}

\begin{theorem}\label{thm:even-conjugacy for twisted covers} Let $r =2t$ be an even integer and let $\rho$ be any lift of $\theta$ to  $\widetilde{G}^{(m)}$. For $\tilde{g}\in \widetilde{G}^{(m)}_{reg}$, there exists $\tilde{y}\in \widetilde{G}^{(m)}$ such that 
\[\rho(\tilde{g})=\tilde{y}\tilde{g}\tilde{y}^{-1}.
\]
\end{theorem}

\begin{theorem}[Dualizing lifts] 
Let $r =2t$ be an even integer and let $G=\GL(r,F)$. Let $\widetilde{G}^{(m)}$ be the $2$-fold Kazhdan-Patterson cover of $G$ with respect to $\Lambda.$
Let $\pi$ be an irreducible admissible genuine representation of $\widetilde{G}^{(m)}$. For any $\alpha \in F^{\times}$, $\sigma_\alpha$ is a dualizing involution, i.e.,
\[\pi^{\sigma_{\alpha}}\simeq \pi^{\vee}. \] 
\end{theorem}

\begin{proof}
 Using Theorem~\ref{thm:even-conjugacy for twisted covers} and proceeding as in Corollary~\ref{conjugacy statement wrt sigma alpha}, it follows that, for
 $\tilde{g} \in \widetilde{G}^{(m)}_{reg}$, there exists $\tilde{u}\in \widetilde{G}^{(m)}$ such that 
\[\sigma_\alpha(\tilde{g})=\tilde{u}\tilde{g}\tilde{u}^{-1}.
\]
Hence, repeating the arguments in Theorem~\ref{main theorem 2} with $\rho$ replaced by $\sigma_\alpha$, it follows that $\sigma_\alpha$ is a dualizing involution. \qedhere
\end{proof}

Recall that $z_{\lambda} = \lambda 1_{G}$ for $\lambda \in F^{\times}$. When $r = 2t+1$ is an odd integer, $pr^{-1}(Z^{(m)})$ is the center of $\widetilde{G}^{(m)}$, where $Z^{(m)}=\{z_{\lambda} : \lambda \in F^{\times} \}$ (see Proposition 0.I.1 in \cite{KazPat[2]}). Thus, given any lift $\rho$  of  $\theta$ to  $\widetilde{G}^{(m)}$, following Section~\ref{Odd rank case}, we can construct an involution $\rho^{\#}$ on $\widetilde{G}^{(m)}$ which satisfies analogous properties as in Lemma~\ref{definition of rhohash} and Lemma~\ref{rho hash fixes zlambda}. 

\begin{theorem}\label{thm:odd-conjugacy for twisted covers}Let $r =2t+1$ be an odd integer. For $\tilde{g}\in \widetilde{G}^{(m)}_{reg}$, there exists $\tilde{y}\in \widetilde{G}^{(m)}$ such that 
\[\rho^{\#}(\tilde{g})=\tilde{y}\tilde{g}\tilde{y}^{-1}.
\]
\end{theorem}

The proof works as in Theorem~\ref{thm:odd-conjugacy}.

\begin{theorem}[Uniqueness] Let $r =2t+1$ be an odd integer and let $G=\GL(r,F)$. Let $\widetilde{G}^{(m)}$ be the $2$-fold Kazhdan-Patterson cover of $G$ with respect to $\Lambda.$ Let $\pi$ be an irreducible admissible genuine representation of $\widetilde{G}^{(m)}$. Then there exists a unique lift $\sigma$ of $\tau$ to $\widetilde{G}^{(m)}$ such that
\[\pi^{\sigma}\simeq \pi^{\vee}.\]
\end{theorem}

\begin{proof}
Since $\tilde{z}_\lambda \in Z(\widetilde{G}^{(m)})$ for every $ \lambda \in F^\times,$ analogue of Lemma~\ref{chirho is trivial iffn} holds. Using Theorem~\ref{thm:odd-conjugacy for twisted covers} and proceeding as in Theorem~\ref{Uniqueness of dualizing involution}, we can show that there exists a unique dualizing lift $\rho$ of $\theta$ to $\widetilde{G}^{(m)}$.  Let $\sigma: \widetilde{G}^{(m)} \to \widetilde{G}^{(m)}$ be the map given by 
\[\sigma(g)=(w_0,1)\rho(g)(w_0,1)^{-1}.\] By similar arguments as in Corollary~\ref{Uniqueness of lift sigma}, it follows that  $\sigma$ is the unique lift of $\tau$ to $\widetilde{G}^{(m)}$ which is dualizing. Hence the result. 
\end{proof}

\def\cprime{$'$}
\providecommand{\bysame}{\leavevmode\hbox to3em{\hrulefill}\thinspace}
\providecommand{\MR}{\relax\ifhmode\unskip\space\fi MR }
\providecommand{\MRhref}[2]{%
  \href{http://www.ams.org/mathscinet-getitem?mr=#1}{#2}
}
\providecommand{\href}[2]{#2}

\end{document}